\documentclass[a4paper,nonatbib]{article}
\pdfoutput=1
\usepackage[nonatbib]{customstyle}
\usepackage{enumitem,graphicx,xurl,xcolor}
\usepackage{booktabs,nicefrac,microtype,wrapfig,multirow,xspace}
\usepackage[titletoc]{appendix}
\def\bi{\begin{itemize}[leftmargin=2em]}
\def\ei{\end{itemize}}
\def\bnum{\begin{enumerate}[leftmargin=2em,label={(\arabic*)}]}
\def\enum{\end{enumerate}}
\def\bmat{\begin{pmatrix}}
\def\emat{\end{pmatrix}}
\def\EE{\mathbb{E}}
\def\CC{\mathbb{C}}
\def\RR{\mathbb{R}}
\def\NN{\mathbb{N}}
\def\PP{\mathbb{P}}
\def\LL{\mathcal{L}}
\def\GG{\mathcal{G}}

\DeclareMathOperator{\rme}{e}
\def\1{\mathbb{I}}
\def\CN{\mathcal{N}_{\CC}}
\def\Var{\mathsf{Var}}
\def\mh{\hspace*{3pt}}
\newcommand{\al}[1]{%
  \begin{align}
    #1
  \end{align}
}
\def\mkakko#1{\left(#1\right)}

\def\kkakko#1{\left[#1\right]}
\def\mmkakko#1{\big(#1\big)}
\renewcommand{\bar}[1]{\overline{#1}}

\renewcommand{\epsilon}{\varepsilon}
\newcommand{\s}[1]{{-#1}}
\let\Re\relax
\let\Im\relax
\DeclareMathOperator{\Re}{Re}
\DeclareMathOperator{\Im}{Im}
\DeclareMathOperator{\Tr}{Tr}
\newcommand{\MP}{Marchenko--Pastur\xspace}
\theoremstyle{lemma}
\newtheorem{lemma}{Lemma}
\theoremstyle{proposition}
\newtheorem{proposition}{Proposition}
\allowdisplaybreaks

\newcommand{\minibox}[1]{%
  \text{%
    \setlength{\fboxsep}{1.3pt}
    \setlength{\fboxrule}{0.4pt}
    \fbox{$\small #1$}
  }%
}

\title{Complex Product Random Matrices with Conjugate Mixing: From Quartic to Sextic Spectral Laws}

\author{%
  Takuya Kanazawa\\
  Department of Business Administration, Kobe Gakuin University\\
  1-1-3 Minatojima, Chuo-ku, Kobe 650-8586, Japan
}

\begin{document}

\maketitle
\thispagestyle{plain}

\begin{abstract}
We study the complex random matrix ensemble 
\[ 
	Q = \frac{1}{\sqrt{d}}
	\bigl(c_1XW+c_2\overline{XW}\bigr)+Z\,, 
\]
where $c_1,c_2\in\mathbb{C}$ are constants and $X\in\mathbb{C}^{n\times d}$, $W\in\mathbb{C}^{d\times p}$ and $Z\in\mathbb{C}^{n\times p}$ are mutually independent random matrices, each with i.i.d.~complex Gaussian entries. The entries may have nonzero pseudo-variances, encoding the second-order couplings between the matrices and their complex conjugates. When $c_2=0$ and all matrices are real, this ensemble reduces to the well-studied product-plus-noise model, which arises in the analysis of shallow neural networks. Our ensemble is a widely linear extension of this single-channel model, incorporating the conjugate channel $\overline{XW}$. Under a standard self-averaging assumption for the resolvents, we perform a leave-one-out analysis to derive self-consistent equations for the asymptotic Stieltjes transform of $Q^*Q/n$ in the proportional limit $n,d,p\to\infty$. Their scalar reduction yields a sextic equation, in contrast to the quartic equation governing the single-channel model. We clarify how this equation reduces to lower-degree equations and recover the Marchenko--Pastur, product-Wishart and Pennington--Worah distributions as special cases. Numerical experiments confirm the theoretical predictions.
\end{abstract}

\clearpage
\hrule\vspace{-0.2em}
\tableofcontents
\vspace{1.6em}
\hrule
\vspace{1.6em}

\section{Introduction}

Products of random matrices arise in a broad range of high-dimensional problems in physics, mathematics, and information science \autocite{Crisanti1993book,Akemann:2011csh,Couillet_Liao_2022}. Examples include multi-stage scattering in multiple-input multiple-output wireless communication \autocite{Muller2002,TulinoBook2004,CDbook}, time-lagged cross-correlations in financial time series \autocite{Bouchaud2007} and quantum chromodynamics at nonzero chemical potential \autocite{Osborn:2004rf,Akemann:2007rf,Kanazawabook}. These applications have motivated extensive studies of both the eigenvalues and singular values of products of random matrices \autocite{Burda2010square,Burda2010,AkemannBurda2012,Akemann2013,AkemannKieburgWei2013,IpsenKieburg2014,Kuijlaars2014,
Dupic2014,Forrester2014product,AkemannIpsen2015review}. Product ensembles display spectral behavior that is substantially richer than that of a single random matrix, since multiplication creates strong dependencies among the entries and may destroy symmetries possessed by the individual factors.

A related class of product-matrix models has recently become important in the analysis of high-dimensional machine learning. Consider a random feature matrix $f(H)$ obtained by applying a nonlinear activation function $f:\RR\to \RR$ entrywise to a preactivation matrix
\al{
	H = \frac{XW}{\sqrt d}\,,
}
with independent real Gaussian matrices $X$ and $W$. Here, $X$ represents the data and $W$ the random weights of a neural network with a single hidden layer. Since the entries of $H$ are not independent, applying an entrywise nonlinearity to $H$ produces a highly correlated random matrix whose spectrum cannot be obtained from ordinary entrywise universality.

In the real-valued setting, a series of works showed that this nonlinear spectral problem can nevertheless be reduced to a more tractable surrogate problem. \textcite{Pennington2017} derived the limiting spectral distribution and obtained a quartic equation for its Stieltjes transform. Subsequent analyses showed that the same limiting spectrum is generated by a model of the form $XW+Z$, in which the nonlinearity is removed at the cost of introducing an independent Gaussian matrix $Z$ \autocite{Peche2019,Benigni2021,Piccolo2021}. Related deterministic-equivalent and learning-theoretic approaches established analogous reductions for random kernel matrices \autocite{El_Karoui_2010,Cheng_Singer_2013} and ridge regression with random features \autocite{Louart2018,Hastie2022,Mei-Montanari-2022}. Block linearization and Dyson equations provide another framework for constructing deterministic equivalents for random feature models and other correlated random matrices \autocite{Erdos2019,LatourelleVigeantPaquette2026}. By contrast, considerably less is known about the spectra of complex-valued random-feature models.

The complex setting introduces a second-order structure with no counterpart in the real setting. A complex random variable is called \emph{proper} when its pseudo-covariance vanishes and \emph{improper} otherwise. For improper data, covariance alone does not provide a complete second-order description. This distinction underlies the widely linear framework in complex signal processing, in which a signal and its complex conjugate are treated as separate inputs \autocite{Brown1969,Picinbono_Chevalier_1995,Mandic2009book,Schreier2010book}. Widely linear methods are now standard for processing noncircular complex signals arising in communications, acoustics, imaging, and related applications. 

The same distinction becomes essential for nonlinear random features in complex-valued neural networks, whose foundations and applications are reviewed in \textcite{Hirose2012book,Lee2022review}. For a general non-holomorphic activation $f:\CC\to\CC$, both the preactivation and its complex conjugate contribute to the corresponding Gaussian surrogate. In a separate work \autocite{kanazawa_cgep_2026}, we argue that the resulting Gaussian equivalent is not the conventional product-plus-noise model; instead, it is a two-channel ensemble involving both $H$ and $\bar H$. Motivated by this observation, we consider 
\al{
	Q = c_1H + c_2 \bar H + Z \,, 
	\qquad H = \frac{XW}{\sqrt d}\,, 
	\label{eq:gsgWQQ}
} 
where $c_1,c_2\in\CC$ are constants and $X\in\CC^{n\times d}$, $W\in\CC^{d\times p}$ and $Z\in\CC^{n\times p}$ are mutually independent complex Gaussian random matrices with \mbox{i.i.d.} entries within each matrix. Their pseudo-variances may be nonzero, giving rise to second-order couplings between the matrices and their complex conjugates. Although motivated by complex-valued neural networks, this conjugate-mixing ensemble \eqref{eq:gsgWQQ} defines a matrix model of independent interest, lying at the intersection of product random matrices and widely linear signal processing. Our main object is the empirical eigenvalue distribution of $Q^*Q/n$ in the proportional high-dimensional limit. 

Under a standard self-averaging assumption for the resolvents, we perform a leave-one-out analysis to derive a closed system of self-consistent equations for the asymptotic Stieltjes transform of $Q^*Q/n$ and the auxiliary quantities induced by the conjugate coupling. In the generic case, the system reduces to a sextic equation for the Stieltjes transform, in contrast to the quartic equation governing the real single-channel model. This increase in algebraic degree reflects the additional correlations carried by the conjugate channel.

The resulting spectral law unifies several known results. Appropriate parameter choices recover the \MP distribution \autocite{Marchenko1967}, the spectral law for a product-Wishart matrix \autocite{Burda2010,Akemann2013}, and the Pennington--Worah distribution \autocite{Pennington2017} for real nonlinear random features. Further parameter constraints produce systematic reductions to equations of lower degree. We classify these reductions and thereby show how established one-channel laws are embedded in the broader conjugate-mixing ensemble. 

We solve the self-consistent equations numerically and compare the resulting densities with empirical spectra of finite-dimensional random matrices, observing close agreement. Together with the Gaussian-equivalence result in \textcite{kanazawa_cgep_2026}, the present spectral theory provides a random matrix foundation for analyzing nonlinear random features in complex-valued neural networks.

The rest of this paper is organized as follows. Section~\ref{sc:prel} introduces the notation for complex random variables and defines the conjugate-mixing ensemble. Section~\ref{sc:spec} states the main result, constructs a block pencil that linearizes the product structure, and derives the asymptotic resolvent. It then reduces the resulting self-consistent equations to a scalar algebraic equation and examines its special cases. Section~\ref{sc:numeric} presents numerical validation of the theoretical spectral densities. Section~\ref{sc:conc} concludes the paper. Technical details are deferred to the appendices.

\section{Preliminaries}
\label{sc:prel}
\subsection{Notation and Definitions}

For any positive integer $K$, the set of integers from $1$ to $K$ is denoted by $[K]$. Throughout this paper, we use $i,j,\dots$ for indices in $[n]$, $\alpha,\beta,\dots$ for indices in $[p]$, and $\mu,\nu,\dots$ for indices in $[d]$.

The sets of real and complex numbers are denoted by $\RR$ and $\CC$, respectively, while $\RR_+$ denotes the set of positive numbers and $\CC_+$ the open upper half-plane. For $z\in\CC$, its real and imaginary parts are denoted by $\Re(z)$ and $\Im(z)$, respectively, and its complex conjugate by $\bar z$. The $d$-dimensional vector with all entries equal to $1$ is denoted by $1_d$.

For a complex matrix $A\in\CC^{m\times n}$, $\bar A$ denotes its complex conjugate, $A^\top$ its transpose, and $A^*$ its conjugate transpose, so that $A^*=\bar{A^\top}$. The $d$-dimensional identity matrix is denoted by $\1_d$. The resolvent and the Stieltjes transform for an $N\times N$ Hermitian matrix $A$ are given by
\al{
	G_A(z) \equiv \frac{1}{A-z\1_N}\,, \qquad
	g_A(z) \equiv \frac{1}{N}\Tr G_A(z)
	\qquad \text{with}\quad z\in\CC_+\,.
}
The subscript $A$ is dropped whenever there is no fear of confusion.

For a complex random variable $\zeta$ with finite second moment, its variance and pseudo-variance are defined by
\al{
	v \equiv \Var(\zeta)
	= \EE\left[\left|\zeta - \EE[\zeta] \right|^2\right],
	\qquad
	\tau \equiv \EE\left[\left(\zeta - \EE[\zeta]
	\right)^2\right],
}
respectively. These quantities satisfy $|\tau|\leq v$. The random variable $\zeta$ is called \emph{proper} if $\tau=0$ and \emph{improper} otherwise. 
Writing $\zeta = \zeta_{\rm R} + i \zeta_{\rm I}$, we have 
\al{
	\tau = \Var(\zeta_{\rm R}) - \Var(\zeta_{\rm I}) 
	+ 2 i \mathsf{Cov}(\zeta_{\rm R},\zeta_{\rm I})\,.
}
Thus, a nonzero pseudo-variance indicates unequal variances of the real and imaginary parts, a nonzero covariance between them, or both.
Circular symmetry implies propriety, while the converse holds for complex Gaussian variables but not for general complex distributions \autocite{Schreier2010book}. In the boundary case $|\tau|=v$, writing $\tau=v\rme^{2i\theta}$, there exists a real-valued random variable $u$ such that $\zeta - \EE[\zeta] = \rme^{i\theta}u$, and hence the distribution of $\zeta$ is supported on a straight line through $\EE[\zeta]$ in the complex plane.

A complex-valued Gaussian distribution with mean $\mu\in\CC$, variance $v\geq 0$, and pseudo-variance $\tau\in\CC$ satisfying $|\tau|\leq v$ is denoted by $\CN(\mu,v,\tau)$. Convergence in probability is denoted by $\xrightarrow{~\PP~}$.

\subsection{Problem Setup}\label{sc:probset}

Let $X\in\CC^{n\times d}$, $W\in\CC^{d\times p}$, and $Z\in\CC^{n\times p}$ be mutually independent random matrices with \mbox{i.i.d.} entries distributed as
\al{
	X_{i\mu}\sim\CN(0,v_x,\tau_x), \qquad 
	W_{\mu\alpha}\sim\CN(0,v_w,\tau_w), \qquad 
	Z_{i\alpha}\sim\CN(0,v_z,\tau_z),
}
respectively. We define the conjugate-mixing ensemble as
\al{
	Q\equiv c_1H + c_2\bar H + Z  \qquad \text{with}\qquad 
	H\equiv\frac{XW}{\sqrt d}\in\CC^{n\times p}\,, \quad 
	c_1,c_2\in\CC\,.
	\label{eq:ensemble}
}

Our primary object of interest is the empirical eigenvalue distribution
of the $p\times p$ matrix $Q^*Q/n$. We study its limiting spectral distribution in the proportional asymptotic regime
\al{
	n,d,p \to\infty \qquad\text{with}\qquad 
	\frac{n}{d}\to \phi\in(0,\infty)\,, \quad
	\frac{p}{d}\to \psi\in(0,\infty)\,.
	\label{eq:lmlm}
}

As a slight generalization, one may add a deterministic constant matrix $\mu 1_n 1_p^\top$ to \eqref{eq:ensemble}. The resulting change in $Q^*Q$ has rank at most two and therefore does not affect the limiting empirical spectral distribution; see \textcite[Theorem~A.44]{Bai2010book} and \textcite[Lemma~2.9]{Couillet_Liao_2022}. It may, however, affect extremal eigenvalues and the associated eigenvectors.

\section{Asymptotic Spectral Analysis}
\label{sc:spec}

\subsection{Main Self-Consistent Equations}\label{sc:mainres}

We now analyze the asymptotic spectrum of the conjugate-mixing ensemble \eqref{eq:ensemble}. The leave-one-out derivation relies on the following assumption.
\paragraph{Assumption.}{\it 
For each fixed $z\in\CC_+$, the normalized traces and matrix-valued partial traces of the resolvents appearing in the derivation converge in probability to nonrandom limits. We further assume that removing any fixed number of rows and columns does not change these limits and that all resolvents involved, including the resulting minor resolvents, are uniformly bounded in operator norm, independently of the matrix dimensions.}\vspace{2mm}

The convergence to deterministic limits reflects the usual self-averaging of macroscopic spectral observables. The remaining conditions ensure that this behavior is stable under the finite-dimensional deletions used in the leave-one-out analysis. Analogous self-averaging and stability properties hold for broad classes of random matrix models \autocite{Mehtabook,Bai2010book,Couillet_Liao_2022,Akemann:2011csh}. Later, in Figure~\ref{fg:mc1}, we provide numerical support for the self-averaging aspect of this assumption by showing that the spectrum of a single large realization agrees remarkably well with the theoretical prediction. 

Throughout this section, all limits of random quantities are understood as limits in probability in the proportional asymptotic regime \eqref{eq:lmlm}. Under the assumption above, we define
\al{
	g_1(z)
	& \equiv \lim_{n,d,p\to\infty} \frac{1}{p}\Tr
	\mkakko{\frac{1}{n}Q^*Q-z\1_p}^{-1},
	\label{eq:g1def}
	\\
	g_2(z) &\equiv \lim_{n,d,p\to\infty}
	\frac{z}{n}\Tr \mkakko{
	\frac{1}{n}QQ^*-z\1_n
	}^{-1}.
	\label{eq:g2def}
}
Note the presence of the factor $z$ in the definition of $g_2(z)$. Since $Q^*Q$ and $QQ^*$ share the same nonzero eigenvalues and differ only in the multiplicities of zero eigenvalues, $g_1(z)$ and $g_2(z)$ satisfy a simple algebraic relation \autocite[Lemma~2.4]{Couillet_Liao_2022}:
\al{
	g_2(z) = \frac{\psi}{\phi} z g_1(z) 
	+ \frac{\psi}{\phi} - 1\,.
	\label{eq:xz0ise0rif}
}

Define the $2\times2$ matrices
\al{
	\Sigma_W = \bmat 
	v_w & \tau_w \\
	\bar{\tau_w} & v_w 
	\emat,
	\qquad \Sigma_{CX} = \bmat
	|c_1|^2v_x & \bar{c_1}c_2\bar{\tau_x}
	\\
	c_1\bar{c_2}\tau_x & |c_2|^2v_x
	\emat.
	\label{eq:sigma_def}
}
Our main result is that, for each fixed $z\in\CC_+$, the following self-consistent equations hold:
\al{
	g_1(z) & = -\frac{1}{
	z + v_zg_2(z) + \Tr\kkakko{\Sigma_W\GG(z)} }\,,
	\label{eq:sce_g1}
	\\
	\GG(z) & = \Gamma_X(z)
	\kkakko{\1_2-\Gamma_W(z)\Gamma_X(z)}^{-1}\,.
	\label{eq:sce_GG}
}
Here,
\al{
	\Gamma_W(z) \equiv \psi g_1(z) \Sigma_W\,,
	\qquad
	\Gamma_X(z) \equiv g_2(z) \Sigma_{CX}\,.
	\label{eq:Gamma_def}
}
Notably, the self-consistent equations depend on the additive noise $Z$ only through its variance $v_z$ and are independent of its pseudo-variance $\tau_z$.

The relevant solution to \eqref{eq:sce_g1} and \eqref{eq:sce_GG} is selected by requiring that $g_1(z)$ be the Stieltjes transform of a probability measure supported on $\RR_{\geq 0}$. Equivalently, it is analytic on $\CC\setminus\RR_{\geq 0}$, maps $\CC_+$ into $\CC_+$, and obeys $zg_1(z)\to-1$ as $|z|\to\infty$.

The definitions of $\Sigma_W$ and $\Sigma_{CX}$ in \eqref{eq:sigma_def} may appear asymmetric, as the coefficients $c_1$ and $c_2$ are incorporated entirely into the $X$ side rather than the $W$ side. To clarify this representation, define $\Sigma_X=\bmat v_x & \tau_x \\ \bar{\tau_x} & v_x \emat$. Then
\al{
	\Sigma_{CX} = \Lambda^*\Sigma_X\Lambda\,,
	\qquad
	\Lambda \equiv \bmat 0 & c_2 \\ c_1 & 0 \emat.
}
An equivalent formulation is obtained by replacing $\{\Sigma_W,\Sigma_{CX}\}$ with $\{\Sigma_{CW},\Sigma_X\}$, where $\Sigma_{CW}\equiv\Lambda\Sigma_W\Lambda^*$. The two formulations lead to the same equation for $g_1(z)$ because the trace term in the denominator of \eqref{eq:sce_g1} is invariant under cyclic permutations.

The derivation of \eqref{eq:sce_g1} and \eqref{eq:sce_GG} is given in the following two subsections. Section~\ref{sc:linearization} constructs a block matrix pencil that linearizes the two-channel matrix products, and Section~\ref{sc:loo} applies the leave-one-out method to derive the self-consistent equations.

\subsection{Linearization and Construction of the Block Pencil}\label{sc:linearization}

\subsubsection{Single-Channel Construction}

To build intuition, we begin by introducing the block pencil in the single-channel case $c_2 = 0$. This preliminary construction illustrates how the target resolvent can be algebraically extracted as a diagonal block of the inverse pencil.

\begin{lemma}
For $z \in \CC_+$, define
\al{
	Q' & \equiv c_1 H + Z = c_1 \frac{XW}{\sqrt{d}} + Z \,,
	\\
	\LL_{\rm base}(z) & \equiv 
	\bmat
		-z \1_p & \frac{1}{\sqrt{n}}Z^* & \frac{1}{\sqrt{d}}W^* & 0 
		\\
		\frac{1}{\sqrt{n}}Z & -\1_n & 0 & \frac{c_1}{\sqrt{n}} X 
		\\
		\frac{1}{\sqrt{d}} W & 0 & 0 & -\1_d 
		\\
		0 & \frac{\bar{c_1}}{\sqrt{n}} X^* & -\1_d & 0
	\emat.
}
Then the upper-left block of $\LL_{\mathrm{base}}(z)^{-1}$ with respect to the displayed $4\times 4$ block partition is
\al{
	\big( \LL_{\rm base}(z)^{-1} \big)_{1,1} 
	= \mkakko{-z\1_p + \frac{1}{n} {Q'}^* Q'}^{-1}.
}
\end{lemma}
\begin{proof}
For brevity, write $\LL_{\text{base}}=\LL_{\text{base}}(z)$. Partition $\LL_{\text{base}}$ into $2 \times 2$ macro-blocks as $\displaystyle \LL_{\mathrm{base}} = \bmat
A & B\\
C & D
\emat$, where
\al{
	A &= 
	\bmat 
		-z \1_p 
		& \frac{1}{\sqrt{n}} Z^* \\ 
		\frac{1}{\sqrt{n}}Z & -\1_n 
	\emat, \quad
	B = 
	\bmat 
		\frac{1}{\sqrt{d}}W^* & 0 
		\\ 
		0 & \frac{c_1}{\sqrt{n}}X 
	\emat, 
	\\
	C &= 
	\bmat 
		\frac{1}{\sqrt{d}}W & 0 
		\\ 
		0 & \frac{\bar{c_1}}{\sqrt{n}} X^* 
	\emat, \quad
	D = 
	\bmat
		0 & -\1_d \\ -\1_d & 0 
	\emat.
}
Taking the Schur complement with respect to $D$ gives 
\al{ 
	\LL_{\text{base}}/D 
	&= A-BD^{-1}C 
	\\ 
	&= \bmat 
		-z\1_p & \frac{1}{\sqrt n}Z^* 
		\\ 
		\frac{1}{\sqrt n}Z & -\1_n 
	\emat 
	+ \frac{1}{\sqrt{nd}} \bmat 
		0 & \bar{c_1} W^*X^* \\ c_1XW & 0 
	\emat 
	= \bmat 
		-z\1_p & \frac{1}{\sqrt n}{Q'}^* 
		\\ 
		\frac{1}{\sqrt n}{Q'} & -\1_n 
	\emat. 
	\label{eq:s1_result} 
} 
By the block inversion formula \eqref{eq:schur_inverse} in Appendix~\ref{app:schur}, the upper-left $2\times2$ macro-block of $\LL_{\text{base}}^{-1}$ is equal to $(\LL_{\text{base}}/D)^{-1}$. Taking the Schur complement of \eqref{eq:s1_result} with respect to its lower-right block $-\1_n$, we obtain 
\al{ 
	\big(\LL_{\text{base}}^{-1}\big)_{1,1} 
	&= \mkakko{ 
		-z\1_p - \frac{1}{\sqrt n} 
		{Q'}^*(-\1_n)^{-1} \frac{1}{\sqrt n}Q' 
	}^{-1} 
	= \mkakko{ -z\1_p+\frac{1}{n} {Q'}^*{Q'} }^{-1}. 
} 
This proves the claimed identity.
\end{proof}

\subsubsection{Conjugate-Mixing Construction}

We now extend the block pencil to the general conjugate-mixing model \eqref{eq:ensemble}. To represent the two conjugate channels in a widely linear form, we introduce the augmented matrices
\al{
	\underline{X} & \equiv \frac{1}{\sqrt{n}}
	\bmat
		c_1 X & c_2 \bar{X}
	\emat
	\in \CC^{n\times 2d}, \qquad 
	\underline{W} \equiv \frac{1}{\sqrt{d}}
	\bmat
		W \\ \bar{W}
	\emat
	\in \CC^{2d\times p}\,. 
	\label{eq:uxuxuwuw}
}
Their product recovers the two product terms in $Q$:
\al{
	\underline{X}\underline{W} &= \frac{1}{\sqrt{nd}}
	\big( c_1 XW + c_2\bar{XW} \big) = \frac{1}{\sqrt{n}}
	\big( c_1 H + c_2 \bar{H} \big) \,.
	\label{eq:43fsnSDX}
}

\begin{lemma}\label{lem:penpen}
For $z\in\CC_+$, define
\al{
	\LL_{\rm full}(z) & \equiv
	\bmat
		-z \1_p & \frac{1}{\sqrt{n}} Z^* & \underline{W}^* & 0
		\\
		\frac{1}{\sqrt{n}} Z & -\1_n & 0 & \underline{X}
		\\
		\underline{W} & 0 & 0 & -\1_{2d}
		\\
		0 & \underline{X}^* & -\1_{2d} & 0
	\emat \,.
	\label{eq:full_pencil_hybrid}
}
Then the $(1,1)$ and $(2,2)$ blocks of $\LL_{\mathrm{full}}(z)^{-1}$ with respect to the displayed $4\times4$ block partition are given by 
\al{ 
	\big( \LL_{\mathrm{full}}(z)^{-1} \big)_{1,1} 
	& = \mkakko{ -z\1_p + \frac1n Q^*Q }^{-1}, 
	\label{eq:lllf1}
	\\
	\big( \LL_{\mathrm{full}}(z)^{-1} \big)_{2,2} 
	& = z\mkakko{ -z\1_n + \frac1nQQ^* }^{-1}. 
	\label{eq:lllf2}
}
\end{lemma}

\begin{proof} 
For brevity, write $\LL_{\mathrm{full}}=\LL_{\mathrm{full}}(z)$ and partition it into $2\times2$ macro-blocks as	$\displaystyle \LL_{\mathrm{full}} = \bmat A & B \\ C & D \emat$, where 
\al{
	A &= \bmat -z\1_p & \frac{1}{\sqrt n}Z^* \\
	\frac{1}{\sqrt n}Z & -\1_n \emat, \qquad 
	B = \bmat \underline W^* & 0 \\ 0 & \underline X \emat, 
	\\
	C &= \bmat \underline W & 0 
	\\
	0 & \underline X^* \emat, \qquad 
	D = \bmat 0 & -\1_{2d} \\ -\1_{2d} & 0 \emat. 
} 
Taking the Schur complement with respect to $D$ gives 
\al{ 
	\LL_{\mathrm{full}}/D &= A-BD^{-1}C 
	\\ 
	&= \bmat -z\1_p & \frac{1}{\sqrt n}Z^* 
	\\ 
	\frac{1}{\sqrt n}Z & -\1_n \emat 
	- \bmat 0 & -\underline W^*\underline X^* 
	\\ -\underline X\underline W & 0 \emat 
	= \bmat -z\1_p & \frac{1}{\sqrt n}Q^* 
	\\ 
	\frac{1}{\sqrt n}Q & -\1_n 
	\emat, 
	\label{eq:LgfLw1wWEF}
} 
where we used \eqref{eq:43fsnSDX}. By the block inversion formula \eqref{eq:schur_inverse} in Appendix~\ref{app:schur}, the upper-left $2\times2$ macro-block of $\LL_{\mathrm{full}}^{-1}$ is the inverse of \eqref{eq:LgfLw1wWEF}. Taking the Schur complement of \eqref{eq:LgfLw1wWEF} with respect to its lower-right block $-\1_n$ gives 
\al{ 
	\big(\LL_{\mathrm{full}}^{-1}\big)_{1,1} 
	&= \mkakko{ 
		-z\1_p - \frac{1}{\sqrt n} 
		Q^*(-\1_n)^{-1} \frac{1}{\sqrt n}Q 
	}^{-1} 
	= \mkakko{ -z\1_p+\frac1nQ^*Q }^{-1}, 
} 
which proves \eqref{eq:lllf1}. Similarly, taking the Schur complement of \eqref{eq:LgfLw1wWEF} with respect to the upper-left block $-z\1_p$ gives 
\al{ 
	\big(\LL_{\mathrm{full}}^{-1}\big)_{2,2} 
	&= \mkakko{ -\1_n - \frac{1}{\sqrt n}Q 
	\mkakko{-z\1_p}^{-1} \frac{1}{\sqrt n}Q^* }^{-1} 
	\\ 
	&= \mkakko{ -\1_n+\frac{1}{nz}QQ^* }^{-1}
	= z\mkakko{ -z\1_n+\frac1nQQ^* }^{-1}, 
} 
which proves \eqref{eq:lllf2}. 
\end{proof}

\subsection{Derivation via the Leave-One-Out Method}
\label{sc:loo}

Under the assumption stated in Section~\ref{sc:mainres}, we derive the self-consistent equations by applying the leave-one-out method to the block pencil. The Schur-complement and leave-one-out identities needed for the derivation are summarized in Appendix~\ref{app:schur}. For a more pedagogical introduction to the leave-one-out method in random matrix theory, we refer to \textcite{Couillet_Liao_2022}.

Let $G(z) = \LL_{\text{full}}(z)^{-1}$ denote the resolvent of the full block pencil $\LL_{\text{full}}(z)$ defined in \eqref{eq:full_pencil_hybrid}. By Lemma~\ref{lem:penpen}, its first two diagonal blocks encode the resolvents of $Q^*Q/n$ and $QQ^*/n$, respectively: 
\al{
	G(z)_{1,1} = \mkakko{ \frac{1}{n}Q^*Q-z\1_p }^{-1}, 
	\quad 
	G(z)_{2,2} = z\mkakko{ \frac{1}{n}QQ^*-z\1_n }^{-1}. 
} 
Recalling \eqref{eq:g1def} and \eqref{eq:g2def}, we have
\al{
	\frac{1}{p}\Tr\mkakko{G(z)_{1,1}} \ 
	\xrightarrow{~\PP~}\ g_1(z)\,, \qquad
	\frac{1}{n}\Tr\mkakko{G(z)_{2,2}} \ 
	\xrightarrow{~\PP~}\ g_2(z)\,.
	\label{eq:gg23874}	
}
For later use, we define the limit of another diagonal block by
\al{
	\GG(z) \equiv \lim_{n,d,p\to\infty}\frac{1}{d} 
	\kkakko{\operatorname{id}_2\otimes\Tr_d}
	\mkakko{G(z)_{3,3}} 
	\,, \qquad \GG(z) \in\CC^{2\times2}\,.
	\label{eq:p0qzpoer}
}
Here, $G(z)_{3,3}$ is regarded as $2\times2$ block matrices whose entries are $d\times d$ matrices, and $\operatorname{id}_2\otimes\Tr_d$ applies the trace to each $d\times d$ block while retaining the outer $2\times2$ structure. In what follows, we write the blocks of $G(z)$ as $G_{11},G_{12},\ldots$ to simplify the notation.

\begin{figure}[tb]
	\centering
	\hspace{13mm}
	\includegraphics[width=5.5cm]{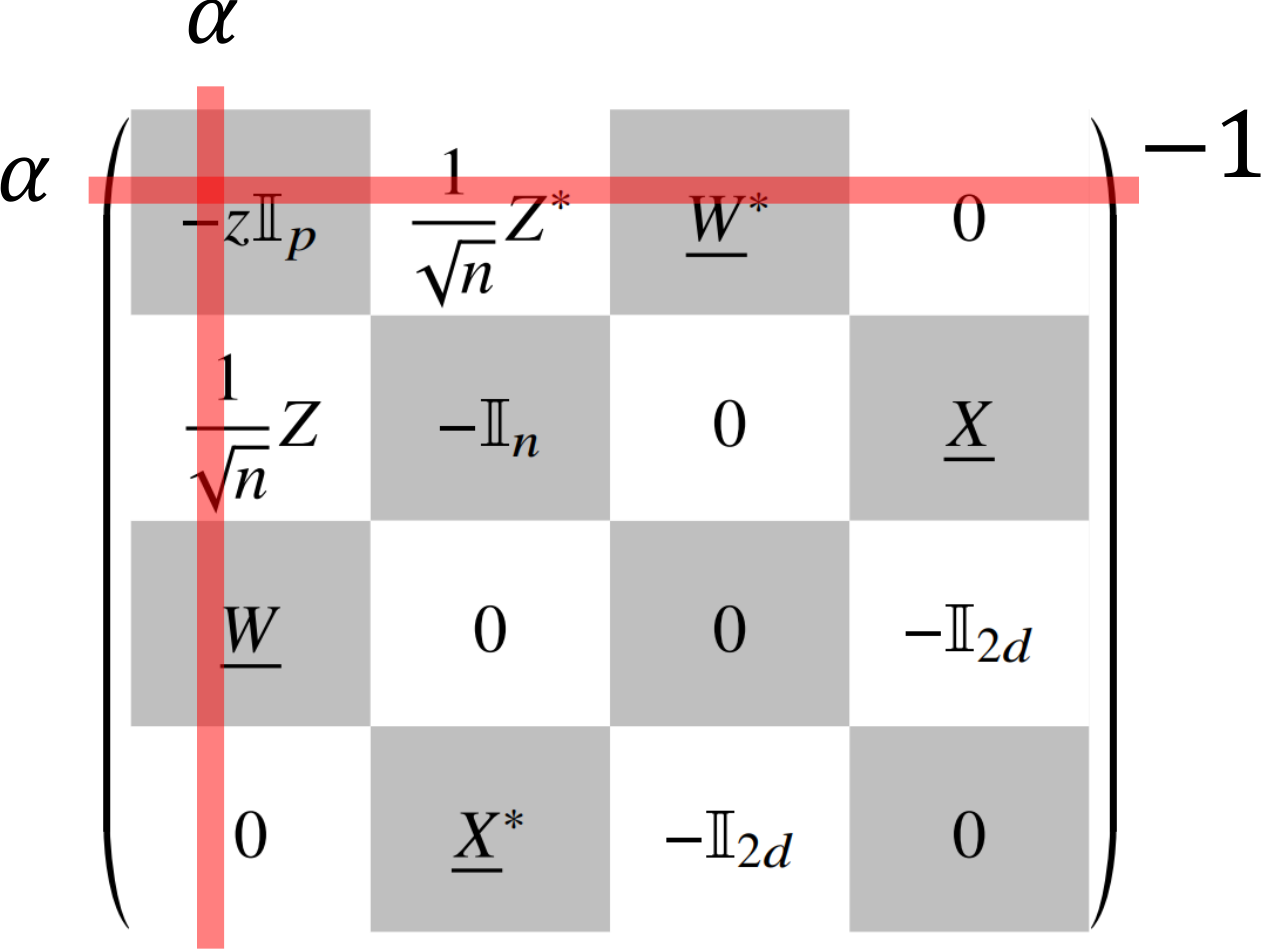}
	\put(-198,48){\Large $G^{\s{\alpha}}=$}
	\caption{\label{fg:loo1}Illustration of the leave-one-out construction for indices in the first block. The red row and column are removed from the underlying matrix $\LL_{\text{full}}(z)$ defined in \eqref{eq:full_pencil_hybrid}.}
\end{figure}

\paragraph{Step 1.} 
We begin by evaluating the diagonal entries of $G_{11}$. For $\alpha\in[p]$, define $G^{\s{\alpha}}$ as the inverse of the submatrix obtained by removing the $\alpha$-th row and column from $\mathcal{L}_{\text{full}}(z)$, as illustrated in Figure~\ref{fg:loo1}.  The blocks of $G^{\s{\alpha}}$ are denoted by $G^{\s{\alpha}}_{11}, G^{\s{\alpha}}_{12}, \dots$ and the sub-blocks inside $G^{\s{\alpha}}_{23},G^{\s{\alpha}}_{33},\dots$ are indicated with the additional indices in $\{1,2\}$; e.g., $\mmkakko{G^{\s{\alpha}}_{33}}_{11}$ refers to the upper-left sub-block of $G^{\s{\alpha}}_{33}$. 

By the leave-one-out formula \eqref{eq:loo_identity_final} in Appendix~\ref{app:schur}, we have, for the scalar diagonal entry $G_{11, \alpha\alpha}$,
\al{
	\mkakko{G_{11,\alpha\alpha}}^{-1} 
	& = -z - \kkakko{
	\bmat
		\frac{Z^*}{\sqrt{n}} & \underline{W}^*
	\emat
	\bmat
		G^{\s{\alpha}}_{22} & G^{\s{\alpha}}_{23}
		\vspace{3pt}\\
		G^{\s{\alpha}}_{32} & G^{\s{\alpha}}_{33}
	\emat
	\bmat
		\frac{Z}{\sqrt{n}} \\ \underline{W}
	\emat
	}_{\alpha\alpha}
	\\
	& = -z - \frac{1}{n} \sum_{i,j=1}^{n}\bar{Z_{i\alpha}}
	G^{\s{\alpha}}_{22,ij}Z_{j\alpha}
	\notag
	\\
	& \quad 
	- \frac{1}{\sqrt{nd}}\sum_{i=1}^{n} \sum_{\mu=1}^{d}
	\Big[
		\bar{Z_{i\alpha}} \mmkakko{G^{\s{\alpha}}_{23}}_{1,i\mu} 
		W_{\mu \alpha}
		+
		\bar{Z_{i\alpha}} \mmkakko{G^{\s{\alpha}}_{23}}_{2,i\mu} 
		\bar{W_{\mu \alpha}}
	\notag
	\\
	& \qquad \qquad 
		+ \bar{W_{\mu\alpha}} \mmkakko{G^{\s{\alpha}}_{32}}_{1,\mu i}
		Z_{i \alpha}
		+ W_{\mu\alpha} \mmkakko{G^{\s{\alpha}}_{32}}_{2,\mu i}
		Z_{i \alpha}
	\Big]
	\notag
	\\
	& \quad 
	- \frac{1}{d}\sum_{\mu,\nu=1}^{d}
	\Big[
		\bar{W_{\mu\alpha}}\mmkakko{G^{\s{\alpha}}_{33}}_{11,\mu\nu}
		W_{\nu\alpha}
		+ \bar{W_{\mu\alpha}}\mmkakko{G^{\s{\alpha}}_{33}}_{12,\mu\nu}
		\bar{W_{\nu\alpha}}
	\notag
	\\
	& \qquad \qquad 
		+ W_{\mu\alpha}\mmkakko{G^{\s{\alpha}}_{33}}_{21,\mu\nu}W_{\nu\alpha}
		+ W_{\mu\alpha}\mmkakko{G^{\s{\alpha}}_{33}}_{22,\mu\nu}
		\bar{W_{\nu\alpha}} 
	\Big]\,.
}

Since $\{Z_{i\alpha}\}_{i\in[n]}$ and $\{W_{\mu\alpha}\}_{\mu\in[d]}$ are independent of $G^{\s{\alpha}}$ by construction, and since all the relevant block matrices have uniformly bounded norms under the assumption of Section~\ref{sc:mainres}, the trace lemma and standard moment estimates for sums of independent centered variables imply, in the limit $n,d,p\to\infty$, that the quadratic terms reduce to their corresponding traces while the cross terms vanish; see \textcite[Lemma~B.26]{Bai2010book} and \autocite[Lemma~2.11]{Couillet_Liao_2022}. For improper complex vectors, we apply the same estimates to their two-dimensional real representations, whose blocks are independent across indices. We thus obtain
\al{
	\mkakko{G_{11,\alpha\alpha}}^{-1} & \xrightarrow{~\PP~} 
	-z -\frac{v_z}{n}\Tr\mmkakko{G^{\s{\alpha}}_{22}}
	- \frac{1}{d} \Big(
		v_w \Tr\kkakko{\mmkakko{G^{\s{\alpha}}_{33}}_{11}}
		+ \bar{\tau_w}\Tr\kkakko{\mmkakko{G^{\s{\alpha}}_{33}}_{12}}
	\notag
	\\
	& \qquad 
	+ \tau_w \Tr\kkakko{\mmkakko{G^{\s{\alpha}}_{33}}_{21}}
	+ v_w \Tr\kkakko{\mmkakko{G^{\s{\alpha}}_{33}}_{22}}
	\Big)
	\\
	& = -z -\frac{v_z}{n}\Tr\mmkakko{G^{\s{\alpha}}_{22}}
	- \frac{1}{d}\Tr\kkakko{
	\bmat 
		v_w & \tau_w \\
		\bar{\tau_w} & v_w 
	\emat 
	\kkakko{\operatorname{id}_2\otimes\Tr_d}
	G^{\s{\alpha}}_{33}} .
}
Under the assumption of Section~\ref{sc:mainres}, the omission of the $\alpha$-th row and column does not alter the limits. Recalling \eqref{eq:gg23874}, \eqref{eq:p0qzpoer} and \eqref{eq:sigma_def}, we find
\al{
	\mkakko{G_{11,\alpha\alpha}}^{-1} & \xrightarrow{~\PP~}
	-z - v_z g_2(z) - \Tr \kkakko{\Sigma_W \GG(z)}.
}
Note that the dependence on $\alpha$ disappears from the right hand side. Averaging the diagonal entries over $\alpha$ and invoking the self-averaging assumption, we obtain
\al{
	g_1(z) & = - \frac{1}{z + v_z g_2(z) + 
	\Tr\kkakko{\Sigma_W \GG(z)}}\,.
	\label{eq:hsAFD6L}
}

\paragraph{Step 2.} 
We next evaluate the diagonal $2\times 2$ block $G_{33}\in\CC^{2d\times2d}$. The standard leave-one-out approach, which removes a single row and column, is unsuitable here. Since the augmented matrices $\underline{W}$ and $\underline{X}$ stack the original variables and their complex conjugates, removing the row and column associated with a single augmented index leaves the corresponding conjugate index in the remaining matrix. This residual correlation prevents a direct application of the trace lemma. 

\begin{figure}[tb]
	\centering
	\hspace{17mm}
	\includegraphics[width=5.5cm]{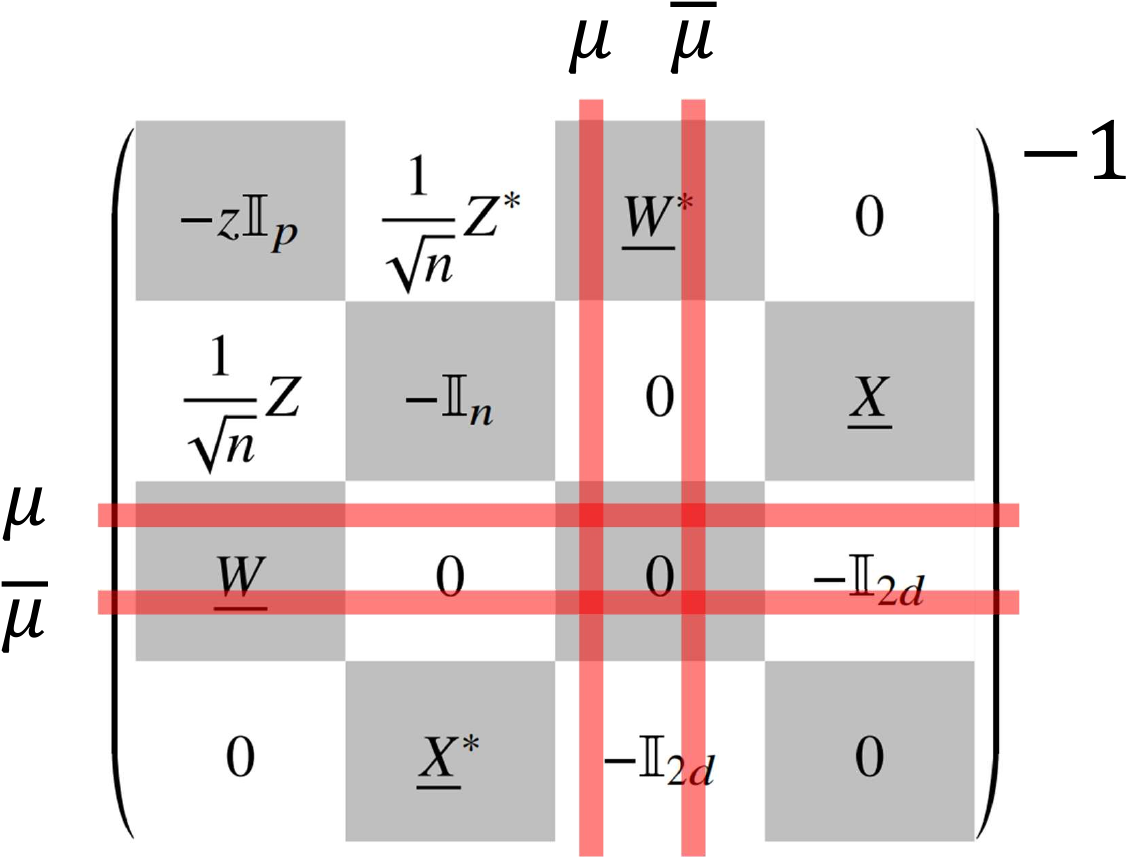}
	\hspace{17mm}
	\includegraphics[width=5.5cm]{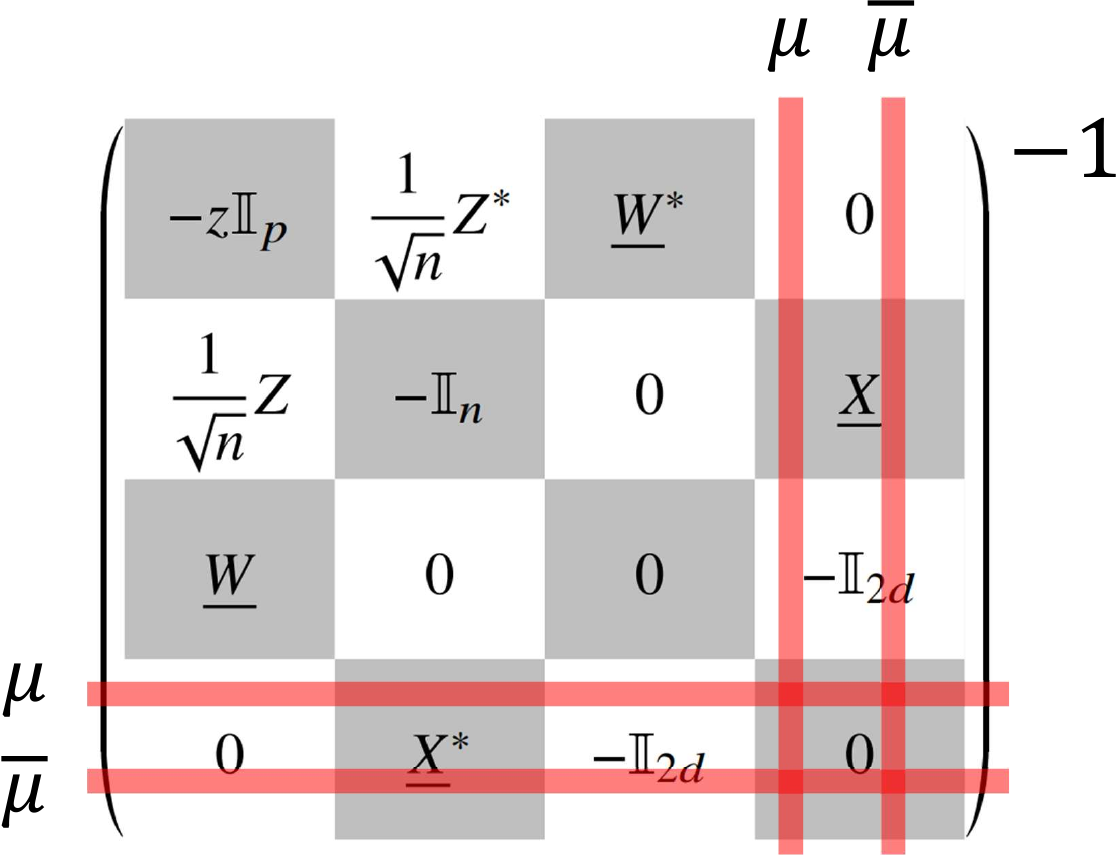}
	\put(-425,47){\Large $G^{\s{(3,\mu)}}=$}
	\put(-210,47){\Large $G^{\s{(4,\mu)}}=$}
	\caption{\label{fg:loo2}Illustration of the leave-one-block-out construction for indices in the third and fourth blocks. The red rows and columns are removed from the underlying matrix $\LL_{\text{full}}(z)$ defined in \eqref{eq:full_pencil_hybrid}.}
\end{figure}

We therefore employ a \emph{leave-one-block-out} approach that simultaneously removes the two rows and columns associated with a variable and its complex conjugate. Denote by $\{\mu,\bar{\mu}\}$ the pair consisting of the $\mu$-th index in the third block of $\mathcal{L}_{\mathrm{full}}(z)$ and its conjugate index. We denote the inverse of the resulting submatrix by $G^{\s{(3,\mu)}}$. Similarly, when the corresponding index pair is removed from the fourth block, we denote the inverse of the resulting submatrix by $G^{\s{(4,\mu)}}$. These constructions are illustrated in Figure~\ref{fg:loo2}. When the corresponding index pairs are removed simultaneously from both the third and fourth blocks, we write $G^{\s{(3,\mu)}\s{(4,\mu)}}$.

We focus on the diagonal $2\times2$ blocks associated with the index pair $\{\mu,\bar{\mu}\}$ in the third and fourth blocks of $\mathcal{L}_{\mathrm{full}}(z)$, denoted by $(G_{33})_{\minibox{\mu}}$ and $(G_{44})_{\minibox{\mu}}$, respectively. Applying \eqref{eq:loo_identity_final} yields
\al{
	\kkakko{(G_{33})_{\minibox{\mu}}}^{-1}
	& = - \kkakko{
	\bmat \underline{W} & -\1_{2d} \emat
	\bmat
		G^{\s{(3,\mu)}}_{11} & G^{\s{(3,\mu)}}_{14}
		\vspace{3pt}
		\\
		G^{\s{(3,\mu)}}_{41} & G^{\s{(3,\mu)}}_{44}
	\emat
	\bmat \underline{W}^* \\ -\1_{2d} \emat
	}_{\minibox{\mu}}
	\\
	& = - \kkakko{
		\underline{W} G_{11}^{\s{(3,\mu)}} \underline{W}^*
		+ G^{\s{(3,\mu)}}_{44}
	}_{\minibox{\mu}} + o^{}_{\PP}(1)
	\\
	& = - \kkakko{
		\frac{1}{d}
		\bmat W \\ \bar{W} \emat G_{11}^{\s{(3,\mu)}} 
		\bmat W^* & W^\top \emat
		+ G^{\s{(3,\mu)}}_{44}
	}_{\minibox{\mu}} + o^{}_{\PP}(1)\,,
}
where $o^{}_{\PP}(1)$ collects the terms that are linear in $\underline{W}$ and vanish in probability in the high-dimensional limit. Using the trace lemma yields 
\al{
	\kkakko{(G_{33})_{\minibox{\mu}}}^{-1}
	& = - \frac{1}{d} \Tr\mkakko{G_{11}^{\s{(3,\mu)}}}
	\Sigma_W - \mkakko{	G^{\s{(3,\mu)}}_{44}
	}_{\minibox{\mu}}
	+ o^{}_{\PP}(1)
	\\
	& = 
	- \psi g_1(z) \Sigma_W - \mkakko{G_{44}^{\s{(3,\mu)}}}
	_{\minibox{\mu}} + o^{}_{\PP}(1)\,.
	\label{eq:34uH0zsdjXX}
}
To evaluate the second term in \eqref{eq:34uH0zsdjXX}, we apply the leave-one-block-out identity once more, now removing the corresponding index pair from the fourth block, and use the trace lemma to obtain
\al{
	\kkakko{\mkakko{G^{\s{(3,\mu)}}_{44}}_{\minibox{\mu}}}^{-1}
	& = - \kkakko{
		\underline{X}^{*} G_{22}^{\s{(3,\mu)}\s{(4,\mu)}}\underline{X}
	}_{\minibox{\mu}}
	\\
	& = - \kkakko{
		\frac{1}{n}\bmat \bar{c_1}X^* \\ \bar{c_2}X^\top \emat
		G_{22}^{\s{(3,\mu)}\s{(4,\mu)}}
		\bmat
			c_1 X & c_2 \bar{X}
		\emat
	}_{\minibox{\mu}}
	\\
	& = - \frac{1}{n}\Tr\mkakko{G_{22}^{\s{(3,\mu)}\s{(4,\mu)}}}
	\Sigma_{CX} + o^{}_{\PP}(1)
	\\
	& \xrightarrow{~\PP~} 
	- g_2(z) \Sigma_{CX} \,.
	\label{eq:Eqwe33034er}
}
Combining \eqref{eq:34uH0zsdjXX} with \eqref{eq:Eqwe33034er} and recalling the definitions of $\Gamma_W(z)$ and $\Gamma_X(z)$ in \eqref{eq:Gamma_def} yields
\al{
	(G_{33})_{\minibox{\mu}}
	& \xrightarrow{~\PP~} \kkakko{- \Gamma_W(z) + \Gamma_X(z)^{-1}}^{-1} 
	\\
	& = \Gamma_X(z)
	\kkakko{\1_2 - \Gamma_W(z)\Gamma_X(z)}^{-1} .
}
When $\Sigma_{CX}$ is singular, namely, when $c_1c_2=0$ or $|\tau_x|=v_x$, the intermediate expression involving $\Gamma_X(z)^{-1}$ is not defined. We interpret it by continuity from the nonsingular case. The final expression above remains well defined because it does not require the inverse of $\Gamma_X(z)$.

Finally, averaging over $\mu$ and recalling \eqref{eq:p0qzpoer}, we obtain
\al{
	\GG(z) = \Gamma_X(z)
	\kkakko{\1_2 - \Gamma_W(z)\Gamma_X(z)}^{-1}.
	\label{eq:Xc1012fD}
}

With \eqref{eq:hsAFD6L} and \eqref{eq:Xc1012fD}, we have completed the derivation of the self-consistent equations.

\subsection{Sextic Equation and Its Reductions}
\label{sc:algebraic_reduction}

In this subsection, we show that the relation \eqref{eq:xz0ise0rif} allows us to eliminate $g_2(z)$ and reduce the matrix-valued self-consistent equations to a single sextic equation for $g_1(z)$. We further classify the principal regimes in which the scalar equation reduces to lower degree, including the quadratic \MP limit.

\begin{proposition}[Scalar sextic equation]
For fixed $z\in\CC_+$, the Stieltjes transform $g_1(z)$ characterized in \eqref{eq:g1def} satisfies the algebraic equation $F(g_1(z))=0$, where
\al{
	F(g) & \equiv T_2 v_z\Xi(g)^3 
	+ \kkakko{T_2\mkakko{\xi(g)+\phi-2}-T_1 v_z} \Xi(g)^2 
	\notag
	\\
	& \qquad + \kkakko{- T_1\mkakko{\xi(g)+\phi-1} +v_z}\Xi(g) 
	+ \xi(g) + \phi
	\label{eq:ZGq7pFD0E}
}
with
\al{
	\xi(g) & \equiv \psi z g + \psi -\phi \,,
	\\
	\Xi(g) & \equiv \frac{\psi}{\phi} g \xi(g)\,,
	\\
	T_1 & \equiv \Tr(\Sigma_W\Sigma_{CX}) 
	= (|c_1|^2 + |c_2|^2) v_xv_w 
	+ 2 \Re(c_1 \bar{c_2}\tau_x\tau_w) \,,
	\label{eq:def_T1}
	\\
	T_2 & \equiv \det(\Sigma_W\Sigma_{CX})
	= |c_1|^2|c_2|^2 (v_x^2 -|\tau_x|^2)
	(v_w^2-|\tau_w|^2) \,.
	\label{eq:def_T2}
}
The equation $F(g)=0$ is generically of degree six in $g$. The branch corresponding to $g_1(z)$ is selected by the Stieltjes conditions in Section~\ref{sc:mainres}.
\end{proposition}

\begin{proof}
For simplicity, we write $g_1=g_1(z)$ and $g_2=g_2(z)$ hereafter. By \eqref{eq:xz0ise0rif}, 
\al{
	g_2 = \frac{\xi(g_1)}{\phi} 
	= \frac{\Xi(g_1)}{\psi g_1}\,. 
	\label{eq:gOT453jef}
}
Combining this expression with \eqref{eq:sce_g1} gives 
\al{
	0 = \frac{1}{g_1}+ z+ v_z \frac{\Xi(g_1)}{\psi g_1}
	+\Tr\kkakko{\Sigma_W\GG(z)}\,. 
	\label{eq:fd843gH} 
} 
For any $2\times2$ matrix $A$, the following identity holds: 
\al{ 
	\Tr\kkakko{ A (\1_2 - xA)^{-1} }
	= \frac{ \Tr A-2x\det A }
	{ 1-x\Tr A+x^2\det A }\,. 
} 
Combining this identity with \eqref{eq:sce_GG}, we obtain 
\al{ 
	\Tr\kkakko{\Sigma_W \GG(z)} &= 
	\Tr\mkakko{ \Sigma_W \Gamma_X(z) 
	\kkakko{\1_2-\Gamma_W(z) \Gamma_X(z)}^{-1} } 
	\\ 
	&= g_2 \Tr\mkakko{ 
	\Sigma_W\Sigma_{CX} 
	\big( \1_2-\psi g_1g_2\Sigma_W
	\Sigma_{CX} \big)^{-1} } 
	\\ 
	&= \frac{\Xi(g_1)}{\psi g_1} 
	\frac{ T_1-2T_2\Xi(g_1) }
	{ 1-T_1\Xi(g_1)+T_2\Xi(g_1)^2 }\,,
	\label{eq:trace_reduction_sextic} 
} 
where $T_1$ and $T_2$ are as defined in \eqref{eq:def_T1} and \eqref{eq:def_T2}, and \eqref{eq:gOT453jef} was used in the last step. Substituting \eqref{eq:trace_reduction_sextic} into \eqref{eq:fd843gH} and multiplying both sides by $\psi g_1$ yield
\al{
	0 & = \psi + \psi zg_1 + v_z \Xi(g_1)
	+ \frac{ \Xi(g_1)\big(T_1-2T_2\Xi(g_1)\big) }
	{ 1-T_1\Xi(g_1)+T_2\Xi(g_1)^2 }
	\\
	& = \phi+\xi(g_1) + v_z \Xi(g_1)
	+ \frac{ \Xi(g_1)\big(T_1-2T_2\Xi(g_1)\big) }
	{ 1-T_1\Xi(g_1)+T_2\Xi(g_1)^2 }\,.
	\label{eq:rt7823Qm4}
}
Finally, clearing the denominator and collecting powers of $\Xi(g_1)$ yield $F(g_1)=0$, where $F$ is given by \eqref{eq:ZGq7pFD0E}.
\end{proof}

\begin{table}
	\caption{Classification of the algebraic reductions of the scalar spectral equation. In the ``Equation'' column, we write $\xi=\xi(g)$ and $\Xi=\Xi(g)$ for brevity. The rightmost column shows a representative matrix form for each case. Here, $X$, $W$ and $Z$ denote independent complex random matrices, while $A$ and $B$ denote independent real random matrices.}
	\label{tb:TT}
	\centering
	\begin{tabular}{|@{}c@{}||@{}c@{}|@{}c@{}|@{}c@{}|@{}c@{}|@{}c@{}|@{}c@{}|@{}c@{}|}
		\hline
		\,Case\, & $T_2$ & $T_1$ & \mh
		$T_1^2-4T_2$\, & $v_z$ 
		& Equation &
		\,$\begin{array}{@{}c@{}}
			\text{Degree}\\[-1pt]\text{in $g$}
		\end{array}$ \hspace{0.5pt}
		& 
		$\begin{array}{c}
			\text{Representative}\\[-1pt]\text{matrix form}
		\end{array}$
		\\[1pt] \hline\hline
		(a) & \multirow{4}{*}[-6mm]{\mh $>0$\mh}
		& \multirow{4}{*}[-6mm]{\mh $>0$\mh}
		& \multirow{2}{*}[-3mm]{$>0$} & \mh$>0$\mh 
		& $F(g)=0$; see \eqref{eq:ZGq7pFD0E} & 6 
		& \mh$a_1 XW + a_2 \bar{XW} + Z$\mh
		\rule{0pt}{12pt}
		\\ \cline{1-1} \cline{5-8}
		(b) && & & $0$ &
		$\begin{array}{c}
			T_2\mkakko{\xi+\phi-2}\Xi^2 
			\\
			- T_1\mkakko{\xi+\phi-1}\Xi
			+ \xi+\phi=0
		\end{array}$
		& 5 & $a_1 XW + a_2 \bar{XW}$
		\rule{0pt}{19pt}
		\\ \cline{1-1} \cline{4-8}
		(c) && & \multirow{2}{*}{$0$} & $>0$ &
		$\begin{array}{c}
			T_1 v_z \Xi^2 
			+ \kkakko{T_1(\xi+\phi-2)-2v_z}\Xi
			\\
			- 2 (\xi + \phi) = 0
		\end{array}$
		& $4$ & $AB + Z$
		\rule{0pt}{19pt}
		\\ \cline{1-1} \cline{5-8}
		(d) &&&& $0$ &
		$\begin{array}{c}
			T_1(\xi+\phi-2)\Xi - 2 (\xi + \phi) = 0
		\end{array}$
		& $3$ & $AB$\rule{0pt}{11pt}
		\\\hline
		(e) &\multirow{4}{*}[0mm]{$0$} 
		& \multirow{2}{*}[0mm]{$>0$} 
		& \multirow{2}{*}[0mm]{$>0$} & $>0$ &
		$\begin{array}{c}
			T_1 v_z \Xi^2 
			+ \kkakko{T_1\mkakko{\xi+\phi-1}-v_z}\Xi
			\\
			-\xi-\phi=0
		\end{array}$
		& 4 & $XW + Z$\rule{0pt}{19pt}
		\\ \cline{1-1} \cline{5-8}
		(f) &&&& $0$ &
		$\begin{array}{c}
			T_1\mkakko{\xi+\phi-1}\Xi -\xi-\phi=0
		\end{array}$
		& 3 & $XW$\rule{0pt}{12pt}
		\\ \cline{1-1} \cline{3-8}
		(g) && \multirow{2}{*}[0mm]{$0$} 
		& \multirow{2}{*}[0mm]{$0$}
		& $>0$ &
		$v_z \Xi + \xi+\phi=0$
		& 2 & $Z$\rule{0pt}{12pt}
		\\ \cline{1-1} \cline{5-8}
		(h) &&&& $0$ & $\xi+\phi=0$ & 1 & 0\rule{0pt}{12pt}
		\\ \hline
	\end{tabular}
\end{table}

Notably, the sextic equation derived above reduces to a lower-degree equation in several special regimes, as summarized in Table~\ref{tb:TT}. To classify these regimes, observe that $\Sigma_W$ and $\Sigma_{CX}$ defined in \eqref{eq:sigma_def} are positive semidefinite. Hence, the matrix $\Sigma_W^{1/2}\Sigma_{CX}\Sigma_W^{1/2}$ has two nonnegative real eigenvalues, denoted by $\lambda_1$ and $\lambda_2$. Since $\Sigma_W\Sigma_{CX}$ has the same eigenvalues, we have 
\al{ 
	T_1 = \lambda_1+\lambda_2 \geq 0\,, \qquad
	T_2 = \lambda_1\lambda_2 \geq 0\,, \qquad 
	T_1^2 - 4T_2 = (\lambda_1-\lambda_2)^2 \geq0\,.
	\label{eq:tT1234r5}
} 
Moreover, $T_2>0$ necessarily implies $T_1>0$. The eight cases in Table~\ref{tb:TT} are classified according to whether $T_1$, $T_2$, $T_1^2-4T_2$, and $v_z$ vanish. Table~\ref{tb:TT} shows that different parameter degeneracies may lead to equations of the same algebraic degree. We discuss these cases below, leaving Cases~(c) and~(d) until the end because they arise on the degeneracy locus $\lambda_1=\lambda_2$.
\bi 
	\item 
	Cases~(a) and~(b) cover the most general setting: widely linear models containing both a product matrix $XW$ and its conjugate $\bar{XW}$, with and without an independent noise matrix $Z$, respectively. Their limiting spectra are governed by sextic and quintic equations. To the best of our knowledge, these spectral laws are new. Since the scalar equations are independent of the pseudo-variance $\tau_z$, they do not distinguish between real and complex $Z$. 
	\item 
	Case~(e) describes a single product channel plus independent noise and defines a broader universality class containing the Pennington--Worah model for real-valued nonlinear random features as a canonical special case \autocite{Pennington2017,Peche2019,Benigni2021}. By \eqref{eq:def_T2}, the condition $T_2=0$ is realized when (i) one of $c_1$ and $c_2$ vanishes; (ii) $X$ becomes real up to a global phase ($|\tau_x|\to v_x$); (iii) $W$ becomes real up to a global phase ($|\tau_w|\to v_w$); or (iv) these mechanisms occur in combination. In each case, the conjugate-mixed response can be reduced to a single effective product. For example, if $X$ is real after absorbing its global phase, then 
	\al{
		c_1XW + c_2\bar{XW} 
		= X \big(c_1W + c_2\bar{W} \big) 
		\equiv XW_{\mathrm{eff}}. 
	} 
	Likewise, if $W$ is real up to a global phase, the response can be written as $X_{\mathrm{eff}}W$. Thus, although the original microscopic models may contain two channels, the condition $T_2=0$ collapses them to a single effective product channel. After matching $T_1$, $v_z$, and the aspect ratios, all these realizations yield the same limiting spectral distribution.
	\item 
	Case~(f), which contains no conjugate-response channel, yields the well-known product-Wishart distribution \autocite{Burda2010,Akemann2013,Dupic2014}. After matching the aspect ratios and the overall variance scale, the cubic equation in Table~\ref{tb:TT} agrees with the corresponding equation in the literature.%
	\footnote{More generally, for a product of $L$ independent random matrices, the Stieltjes transform of the limiting squared singular-value distribution satisfies a polynomial equation of degree $L+1$ \autocite{Burda2010}; here $L=2$.}
	\item 
	Case~(g) yields the \MP distribution for the sample covariance matrix $\frac{1}{n}Z^*Z$ \autocite{Marchenko1967}. 
	\item 
	Case~(h) corresponds to the zero matrix and yields the trivial solution $g(z)=-1/z$.
	\item Finally, we turn to Cases~(c) and~(d), which arise on the degeneracy locus $T_1^2-4T_2=0$. In this regime, the two eigenvalues $\lambda_1$ and $\lambda_2$ coalesce, and the sextic equation \eqref{eq:ZGq7pFD0E} factorizes as 
	\al{
		F(g) = \frac{1}{4}\big(2 - T_1\Xi(g)\big) 
		\mkakko{ -T_1v_z\Xi(g)^2 
		+\kkakko{-T_1(\xi(g)+\phi-2)+2v_z}\Xi(g) 
		+2(\xi(g)+\phi) }\,. 
		\label{eq:AG43DF9} 
	} 
	The roots of $2 - T_1\Xi(g) = 0$ are not admissible because they correspond to poles of \eqref{eq:rt7823Qm4}. In fact, under this condition, \eqref{eq:rt7823Qm4} becomes
	\al{
		0 = \phi+\xi(g) + v_z \Xi(g) + 
		\frac{2 T_1 \Xi(g)}{2-T_1\Xi(g)}\,.
	}
	The valid equation is therefore given by the quartic factor in \eqref{eq:AG43DF9}. Its fourth-order term vanishes in Case~(d), so the resulting degrees in Cases~(c) and (d) are four and three, respectively, as reported in Table~\ref{tb:TT}. In Appendix~\ref{app:fac}, we show that the condition $T_1^2-4T_2=0$, under $T_2>0$, holds if and only if 
	\al{ 
		|c_1|=|c_2|>0 \,,\qquad 
		\frac{|\tau_x|}{v_x} = \frac{|\tau_w|}{v_w} \,, 
		\quad \text{and}
		\quad c_1\bar{c_2}\tau_x\tau_w\in\RR_{\leq 0}\,. 
	} 
	Under these conditions, the widely linear response $c_1XW + c_2\bar{XW}$ may be represented as $AB$, where $A$ and $B$ are independent real random matrices and their inner dimension is $2d$ rather than $d$ (see Appendix~\ref{app:fac}). After matching the overall variance scale, Cases~(c) and~(d) reduce to the single-channel Cases~(e) and~(f), respectively, with the aspect ratios replaced by $(\phi/2,\psi/2)$. Indeed, under this replacement, $\xi$ and $\Xi$ are both divided by two, and the equations in Cases~(e) and~(f) reproduce those in Cases~(c) and~(d). Thus, although Cases~(c) and~(d) arise from the special degeneracy condition $T_1^2-4T_2=0$, they can be absorbed into Cases~(e) and~(f) with doubled inner dimension. By contrast, away from this degeneracy locus, the general conjugate-mixing model retains its two-channel structure and is governed by the sextic or quintic equation in Cases~(a) and~(b).
\ei

\section{Numerical Validation}
\label{sc:numeric}
\subsection{Computational Procedure}

In general, the coupled deterministic equations \eqref{eq:xz0ise0rif}, \eqref{eq:sce_g1} and \eqref{eq:sce_GG} cannot be solved explicitly in a closed algebraic form. Instead, for any given evaluation point $z = x + i\eta$ ($\eta > 0$), the value of $g_1(z)$ can be computed numerically by a damped fixed-point iteration. Here, we detail the algorithmic procedure and the subsequent post-processing steps required to extract the limiting spectral density $\rho(x)$ from the converged Stieltjes transforms.

\paragraph{Step 1. Fixed-Point Iteration}
For each $x\in\RR$, we choose a small $\eta > 0$ and set $z = x + i\eta$. Denoting the iteration step by $t = 0, 1, 2, \dots$, the iterative update proceeds as follows: First, we initialize $g_1^{(0)} = i$ and determine $g_2^{(0)}$ via the algebraic relation \eqref{eq:xz0ise0rif}. Given the tentative solutions $g_1^{(t)}$ and $g_2^{(t)}$ at step $t$, the transition to step $t+1$ is performed through the following steps:
\bnum
	\item Construct the $2 \times 2$ auxiliary matrices:
	\al{
		\Gamma_W^{(t)} = \psi g_1^{(t)} \Sigma_W\,, 
		\quad 
		\Gamma_X^{(t)} = g_2^{(t)} \Sigma_{CX}\,.
	}
	\item Update the $2 \times 2$ resolvent $\GG^{(t)}$ via matrix inversion:
	\al{
		\GG^{(t)} &= \Gamma_X^{(t)} 
		\mkakko{
			\1_2 - \Gamma_W^{(t)} \Gamma_X^{(t)}
		}^{-1}\,,
	}
	\item Compute the new target values. To improve numerical stability, we introduce a damping factor $\delta\in (0, 1)$ to form a convex combination of the previous iterate and the new target value:
	\al{
		g_1^{\text{new}} & = - \frac{1}{z + v_z g_2^{(t)} + 
		\Tr \big[ \Sigma_W \GG^{(t)} \big]} \,,
		\\
		g_1^{(t+1)} & = (1 - \delta) g_1^{(t)} + \delta g_1^{\text{new}}\,,
		\\ 
		g_2^{(t+1)} & = \frac{\psi}{\phi}z g_1^{(t+1)} 
		+ \frac{\psi}{\phi} - 1 \,.
	}
\enum
These updates are repeated until $\max_{i=1,2}|g_i^{(t+1)} - g_i^{(t)}|$ falls below a predefined tolerance. When solving continuously along a finely discretized grid of $x$, we apply a warm start (initializing the iteration at $x_{k+1}$ with the converged solution from the adjacent point $x_k$), which both accelerates convergence and facilitates continuous tracking of the Stieltjes branch.%
\footnote{The relevant branch is selected by the Stieltjes conditions, in particular by mapping $\CC_+$ into $\CC_+$ and satisfying the appropriate behavior at infinity. In our numerical implementation, initialization in the upper half-plane, together with damping and continuation in $x$, reliably tracks this branch.}

\paragraph{Step 2. Extracting the Density via Stieltjes Inversion}
The macroscopic spectral density $\rho(x)$ of $Q^*Q/n$ can be extracted from the converged function $g_1(z)$ by approaching the real axis from the upper half-plane: 
\al{
	\rho(x) = \frac{1}{\pi} \lim_{\eta \to 0^+} \Im[g_1(x + i\eta)]\,. 
}
Since this limit cannot be taken exactly in numerical evaluation, we use a sufficiently small but fixed value of $\eta>0$. 

\paragraph{Step 3. Removal of the Zero-Eigenvalue Contribution}
While setting $\eta>0$ is computationally necessary, it regularizes the spectral measure into a smooth density. This can create a pronounced peak near the origin reflecting the presence of macroscopically many zero eigenvalues of $Q^*Q$. Let $r$ denote the rank of $Q^*Q$. When $v_z>0$, almost surely $r=\min(p,n)$. When $v_z=0$, additional care is needed. If $T_2>0$, the factorization \eqref{eq:43fsnSDX}, together with the nondegeneracy of the two augmented factors, gives $r=\min(p,2d,n)$ almost surely. If $T_2=0$ and $T_1>0$, we instead have $r=\min(p,d,n)$. When $T_1=0$, we have $Q=0$ and hence $r=0$. These rank formulas are also consistent with the zero-mode fractions obtained by successively taking the limits $v_z\to 0$ and $z\to 0$ in the sextic equation \eqref{eq:ZGq7pFD0E}. As a consequence, the $p-r$ trivial zero eigenvalues give rise to a Lorentzian peak centered at the origin, artificially lifting the nonzero spectral density near the origin. Since the empirical spectral measure is normalized by $p$, the atom at zero has mass $(p-r)/p$, yielding
\al{
	\rho_{\text{zero}}(x) = \frac{p - r}{p} 
	\frac{\eta}{\pi (x^2 + \eta^2)}\,.
}
To compare with a histogram of empirically computed nonzero eigenvalues, we subtract this contribution and renormalize the remaining mass to unity, obtaining
\al{
	\rho_{\text{nonzero}}(x) = 
	\max \left( 0,~ \frac{p}{r} \kkakko{
		\frac{1}{\pi}\Im[g_1(x + i\eta)] 
		- \frac{p - r}{p} \frac{\eta}{\pi (x^2 + \eta^2)}
	} \right),
}
where the outer maximum is a numerical clipping used to suppress small negative values caused by finite-$\eta$ and discretization errors.

In numerical implementation of the procedure above, we use $\eta=3\times10^{-5}$, $\delta = 0.5$, and a convergence tolerance of $10^{-7}$ throughout.

\subsection{Comparison with Empirical Spectra}

We now compare the empirical spectra with the corresponding theoretical predictions. All code used to generate the numerical results in this subsection is publicly available on GitHub.%
\footnote{\url{https://github.com/TaTKSM/conjugate-mixing-rmt}} For clarity, we adopt the following two conventions:
\bi
	\item 
	We plot the singular values of $Q/\sqrt{n}$, rather than the eigenvalues of $Q^*Q/n$. The latter are simply the squares of the former. The theoretical density is transformed accordingly.
	\item 
	Assuming $c_1\neq0$, we use the rescaling and phase freedoms of the model to set, without loss of generality,
	\al{
		c_1=v_x=v_w=1\,.
	}
	We further choose the phases of $X$ and $W$ so that $\tau_x$ and $\tau_w$ are real and nonnegative. The remaining free parameters are therefore
\al{
	c_2\in\CC,\qquad 
	(\tau_x,\tau_w)\in[0,1]^2,\qquad
	v_z\geq0,\qquad
	\quad |\tau_z|\leq v_z, \qquad (n,d,p)\in\NN^3.
}
\ei

We display the normalized spectral density averaged over $N_{\rm MC}$ independent Monte Carlo realizations, with $N_{\rm MC}$ indicated in each plot. In most cases, we use $N_{\rm MC}=300$.

\begin{figure}[tb]
	\centering
	\includegraphics[width=.325\textwidth]{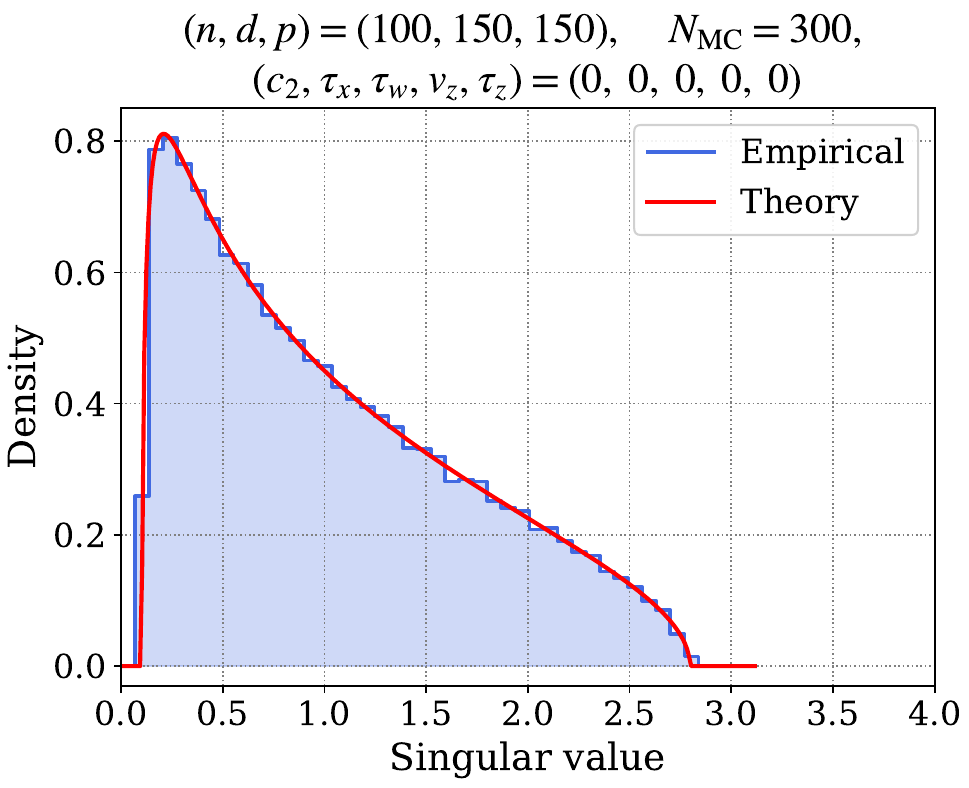}
	\includegraphics[width=.325\textwidth]{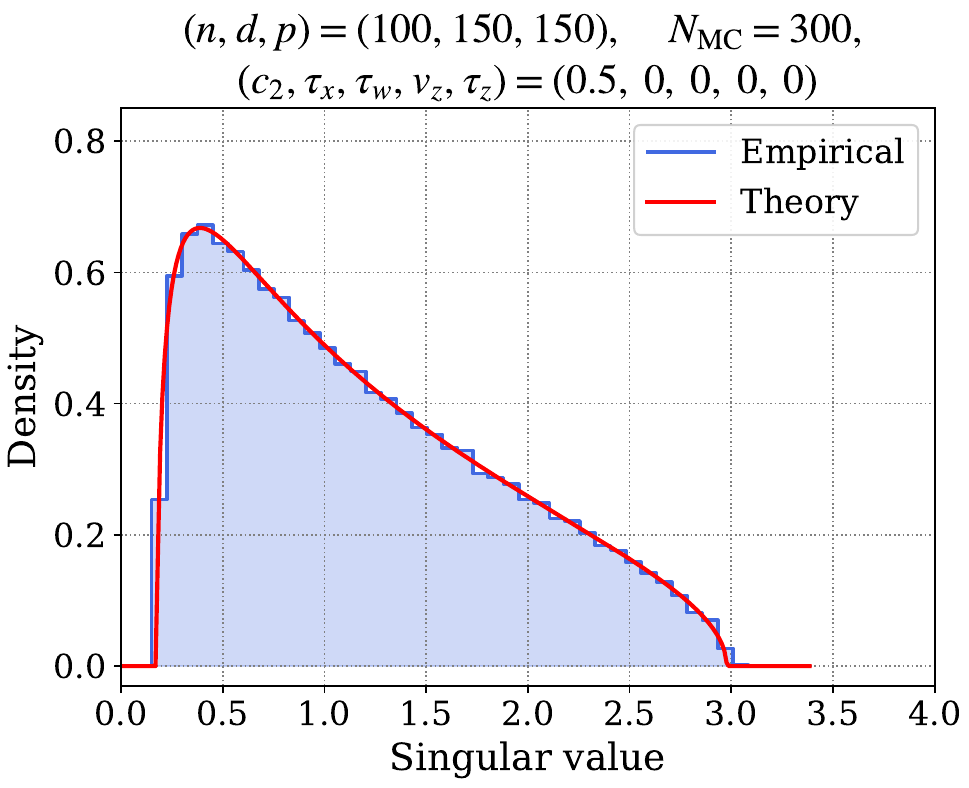}
	\includegraphics[width=.325\textwidth]{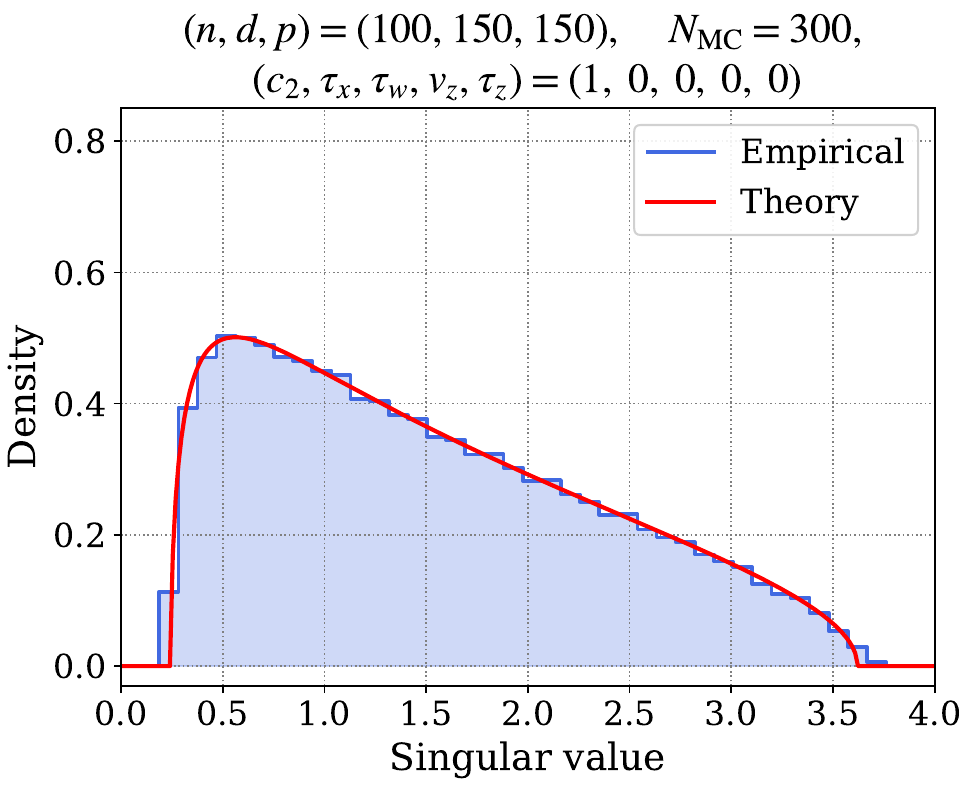}
	\vspace{3mm}\\
	\includegraphics[width=.325\textwidth]{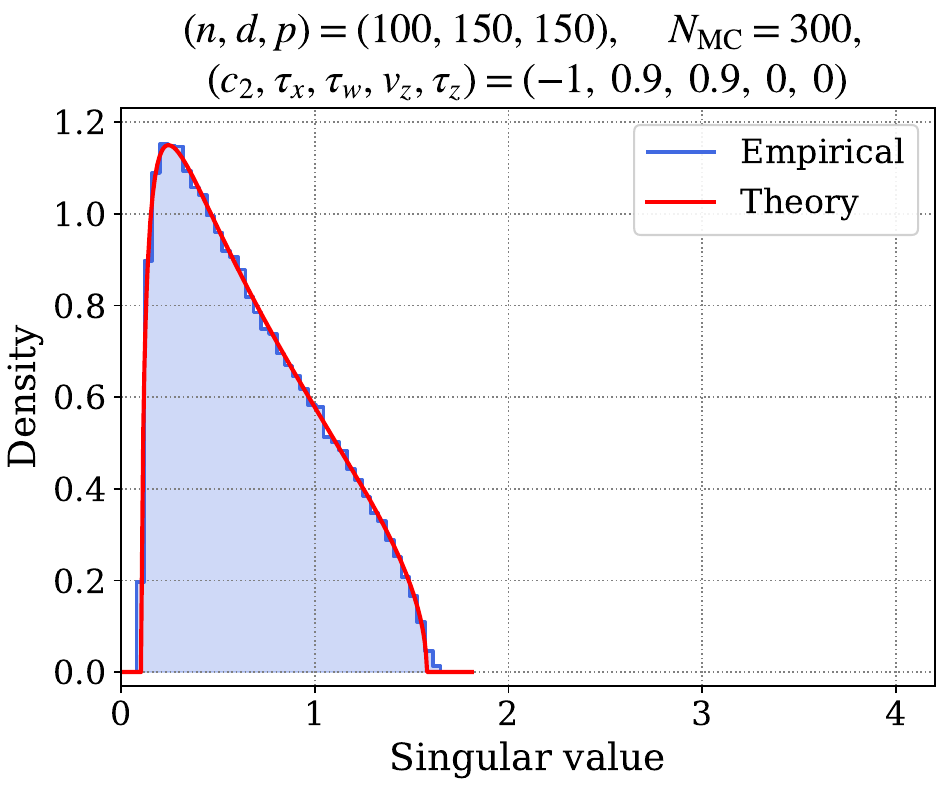}
	\includegraphics[width=.325\textwidth]{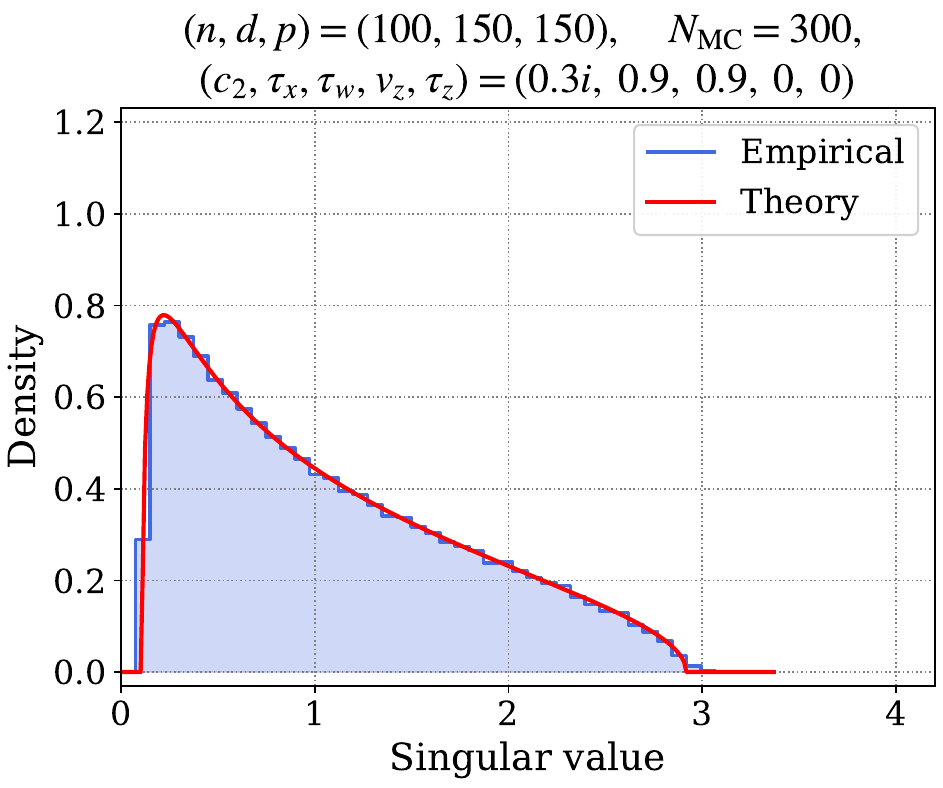}
	\includegraphics[width=.325\textwidth]{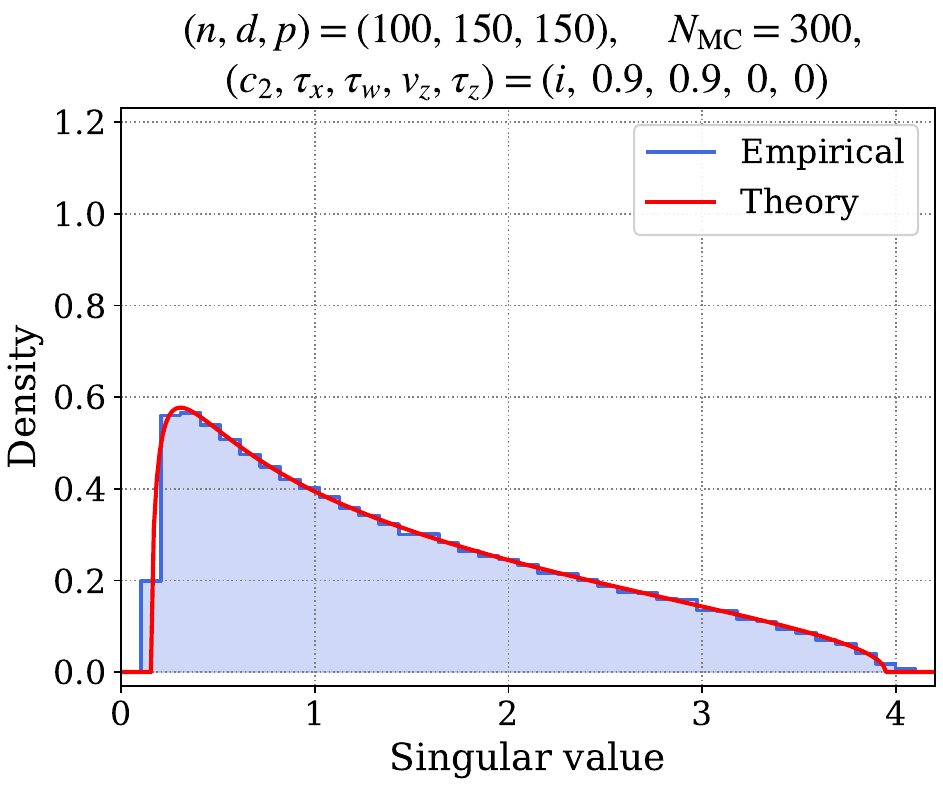}
	\caption{\label{fg:0tr}	Empirical spectra and theoretical predictions for varying $c_2$ without noise ($v_z=0$). Top row: $X$ and $W$ are proper ($\tau_x=\tau_w=0$). Bottom row: $X$ and $W$ are improper ($\tau_x=\tau_w=0.9$).}
\end{figure}

\begin{figure}[tb]
	\centering
	\includegraphics[width=.325\textwidth]{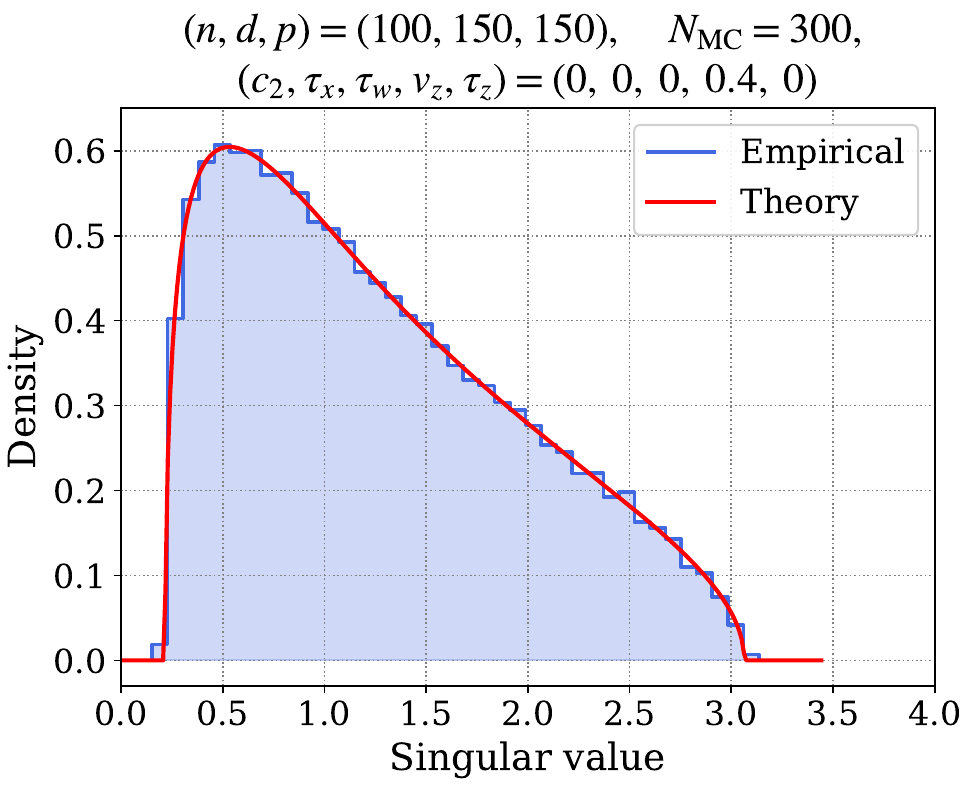}
	\includegraphics[width=.325\textwidth]{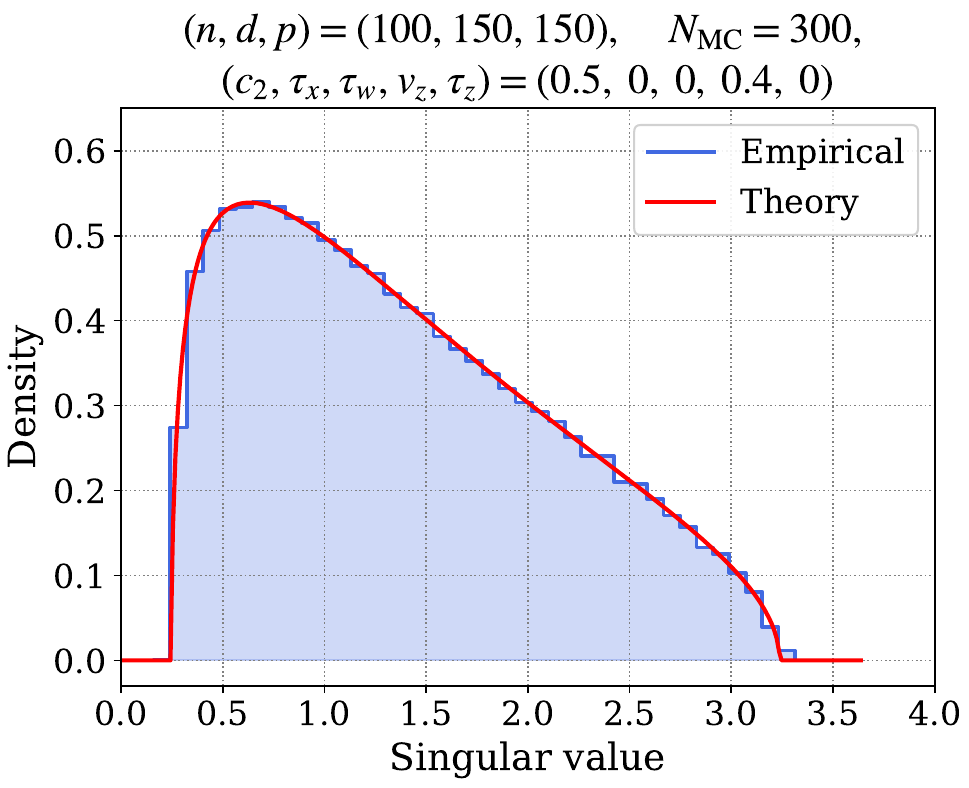}
	\includegraphics[width=.325\textwidth]{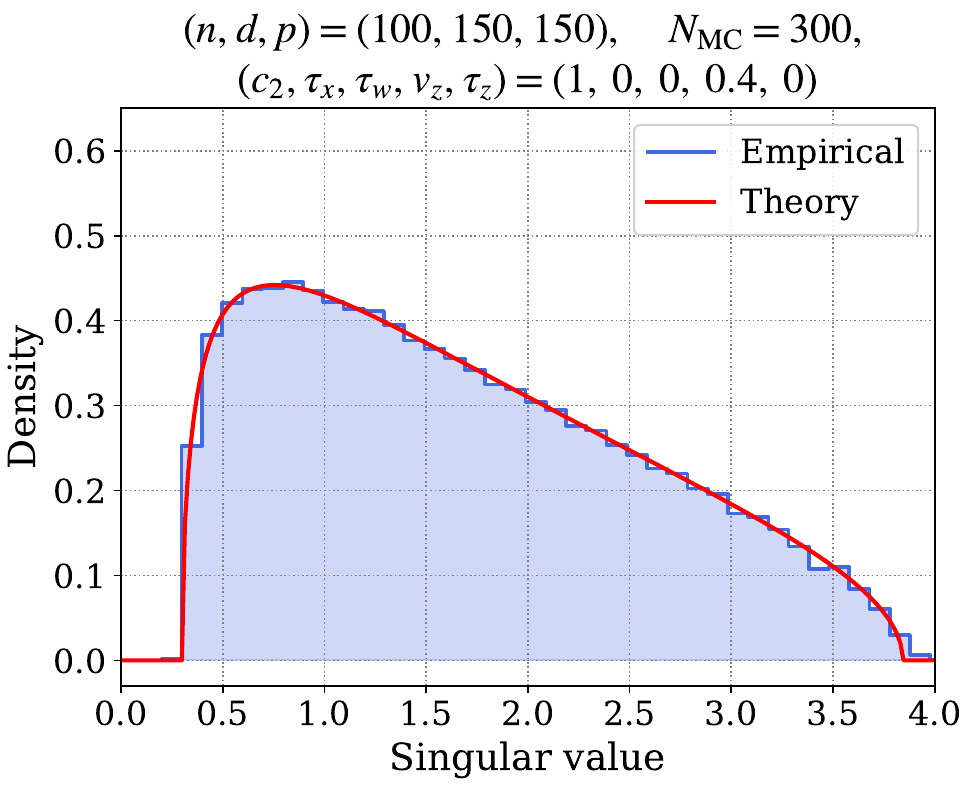}
	\caption{\label{fg:c2vvzz}	Dependence of the spectrum on $c_2$ in the presence of noise ($v_z=0.4$).}
\end{figure}

Figure~\ref{fg:0tr} shows how the spectrum varies with $c_2$ in the absence of noise, $v_z=0$. When $X$ and $W$ are proper (top row), the phase of $c_2$ can be absorbed into the phase freedom of $X$ or $W$, so that $c_2$ may be taken to be nonnegative without loss of generality. As $c_2$ increases, the entrywise variance of $Q$ increases, and the spectral support expands accordingly. When $X$ and $W$ are improper, the spectrum depends more intricately on both the magnitude and phase of $c_2$ (bottom row). In all cases, we observe excellent agreement between our asymptotic theory and the empirical results at finite matrix dimensions. The overall behavior remains qualitatively similar when the noise matrix $Z$ is added, as shown in Figure~\ref{fg:c2vvzz}.

\begin{figure}[tb]
	\centering
	\includegraphics[width=.325\textwidth]{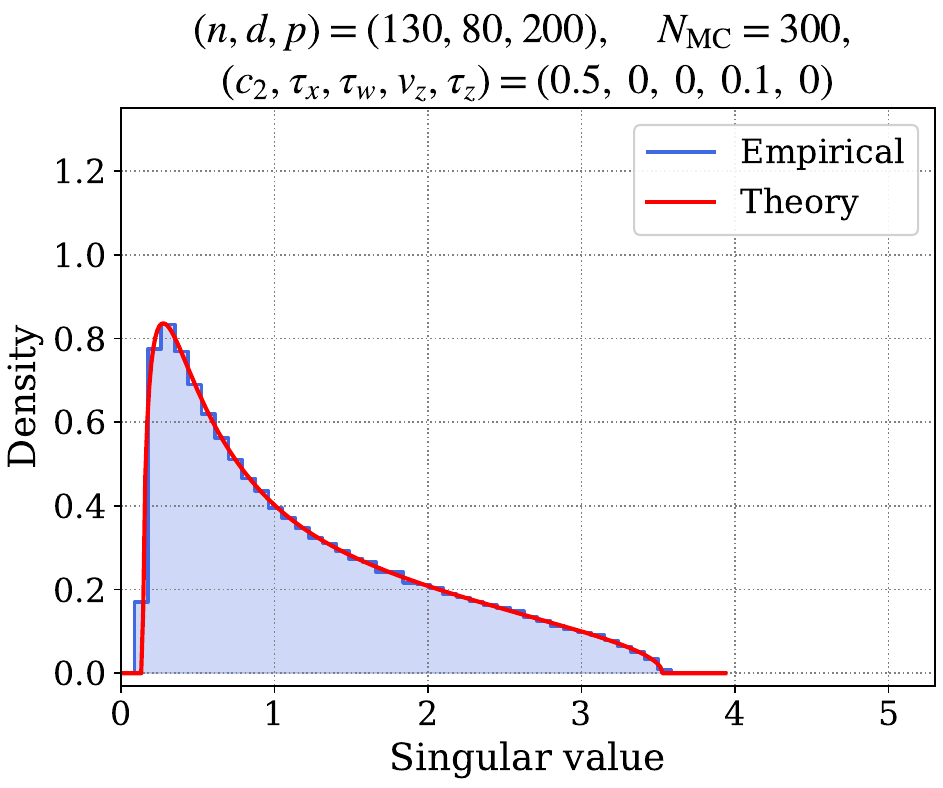}
	\includegraphics[width=.325\textwidth]{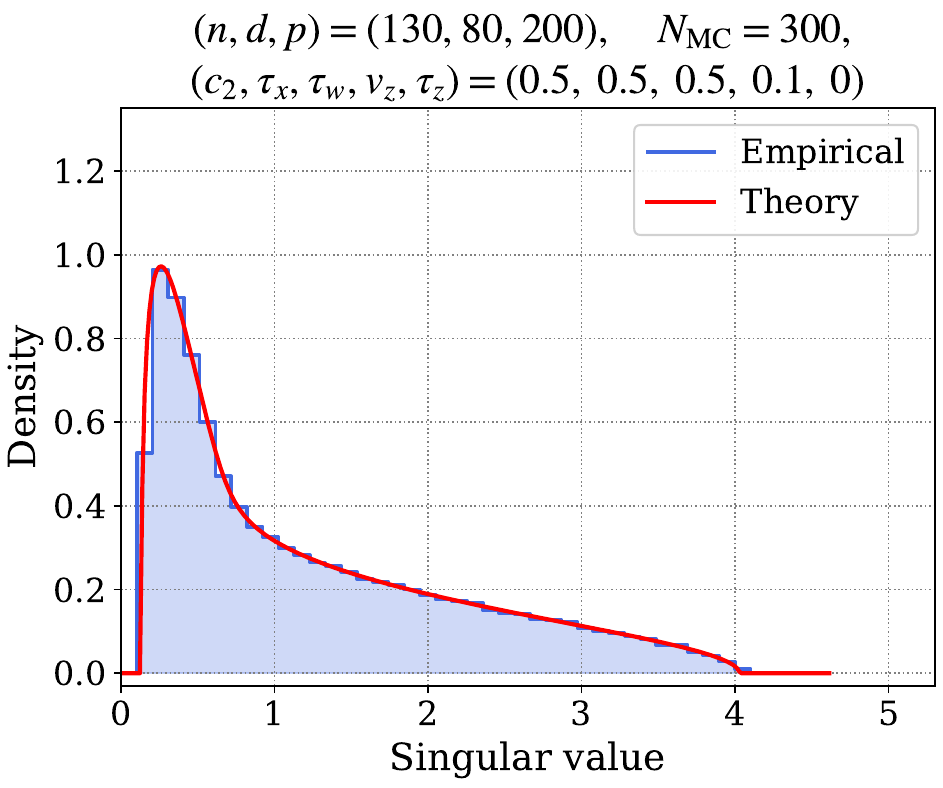}
	\includegraphics[width=.325\textwidth]{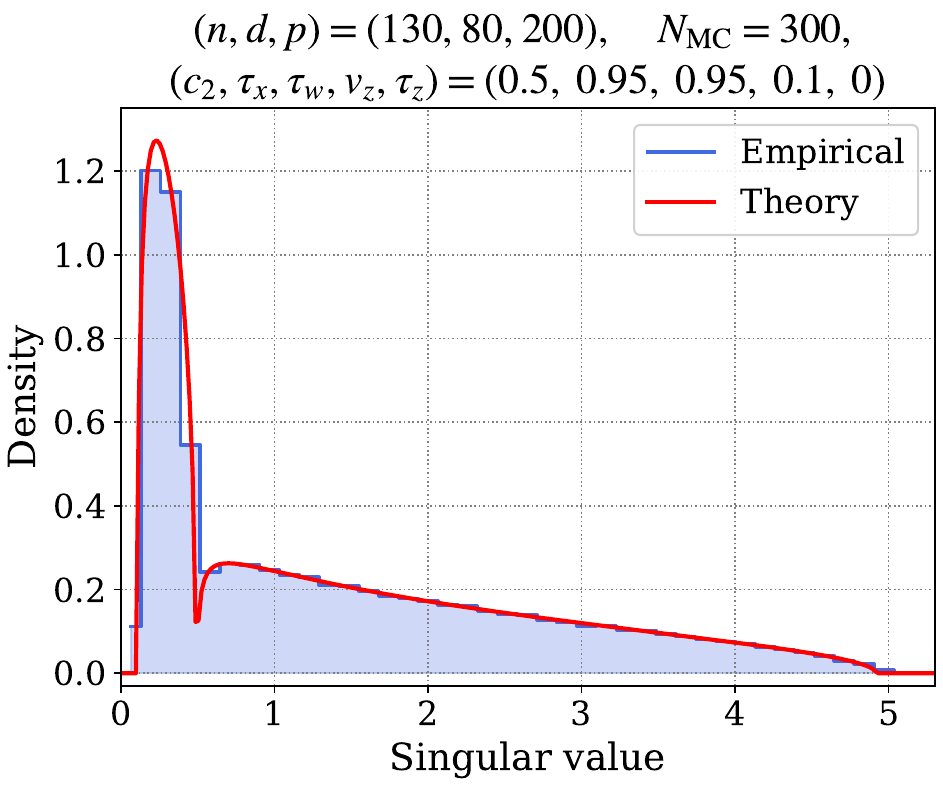}
	\caption{\label{fg:ttt}	Dependence of the spectrum on the impropriety parameters $\tau_x$ and $\tau_w$.}
\end{figure}

Figure~\ref{fg:ttt} shows how the spectrum varies as $\tau_x=\tau_w$ takes the values $0$, $0.5$, and $0.95$. Interestingly, as these parameters approach unity, the density gradually develops a dip and separates into a sharp peak near the origin and a broader peak away from it. This two-bump structure has a simple interpretation in terms of the rank reduction that occurs at the degenerate boundary. For simplicity of exposition, first consider the noise-free case $v_z=0$. By \eqref{eq:uxuxuwuw},
\al{
	Q=\sqrt{n}\,\underline{X}\underline{W}\,,
	\qquad
	\underline{X}\in\CC^{n\times 2d}\,,
	\qquad
	\underline{W}\in\CC^{2d\times p}\,.
}
For $|\tau_x|<v_x=1$ and $|\tau_w|<v_w=1$, the two augmented channels are nondegenerate, and hence
\al{
	\operatorname{rank}(Q)
	= \min(n,2d,p) =\min(130,160,200) = 130
}
almost surely. At the boundary $|\tau_x|=|\tau_w|=1$, however, $X$ and $\overline{X}$ become linearly dependent up to a global phase, and the same occurs for $W$ and $\overline{W}$. The two augmented channels then collapse into a single effective product
channel, giving
\al{
	\operatorname{rank}(Q) = \min(n,d,p)
	= \min(130,80,200) = 80\,.
}
For every $|\tau_x|,|\tau_w|<1$, all $130$ nonzero singular values remain strictly positive. As $\tau_x$ and $\tau_w$ approach unity, however, $50$ of these $130$ singular values become increasingly small and move toward zero, while the remaining $80$ stay away from the origin. This splitting of the singular values naturally produces the two bumps observed in the rightmost panel of Figure~\ref{fg:ttt}. For the small but nonzero value of $v_z$ used in the figure, the noise prevents the smaller singular values from reaching zero exactly, but the same separation into two spectral scales remains clearly visible.

\begin{figure}[tb]
	\centering
	\includegraphics[width=.325\textwidth]{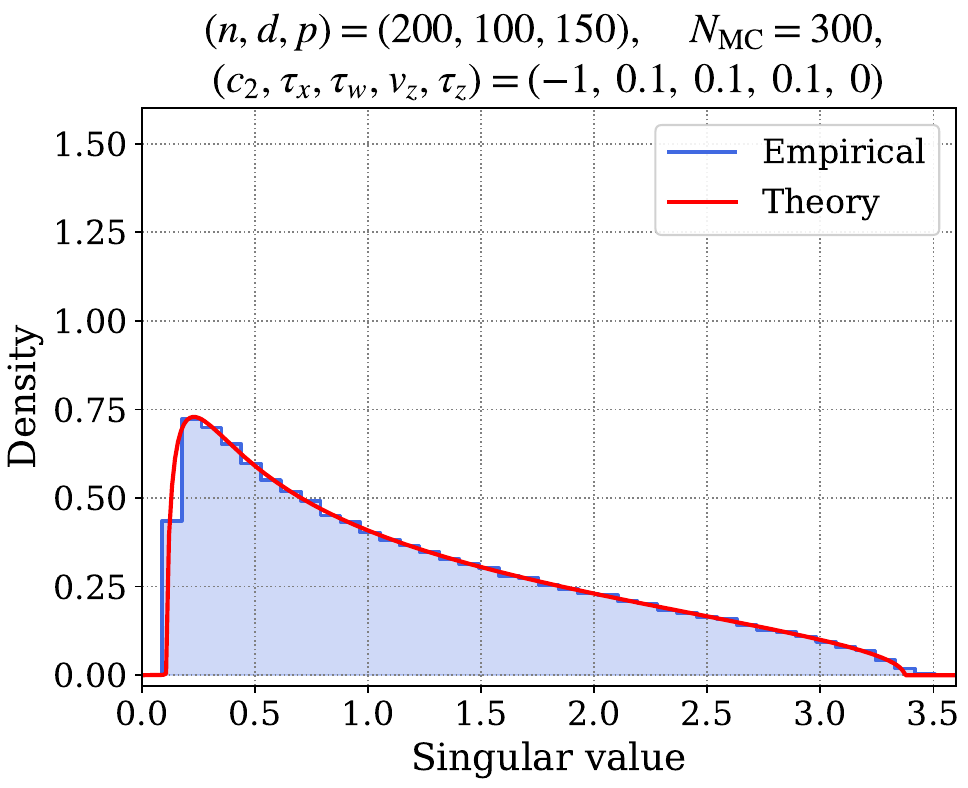}
	\includegraphics[width=.325\textwidth]{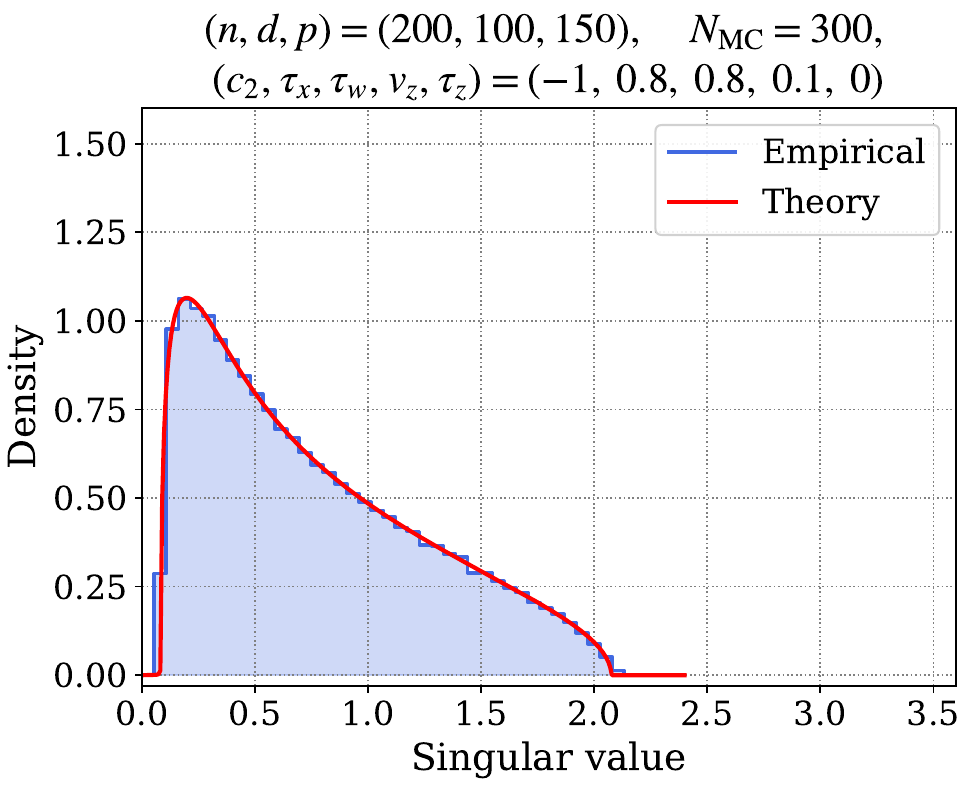}
	\includegraphics[width=.325\textwidth]{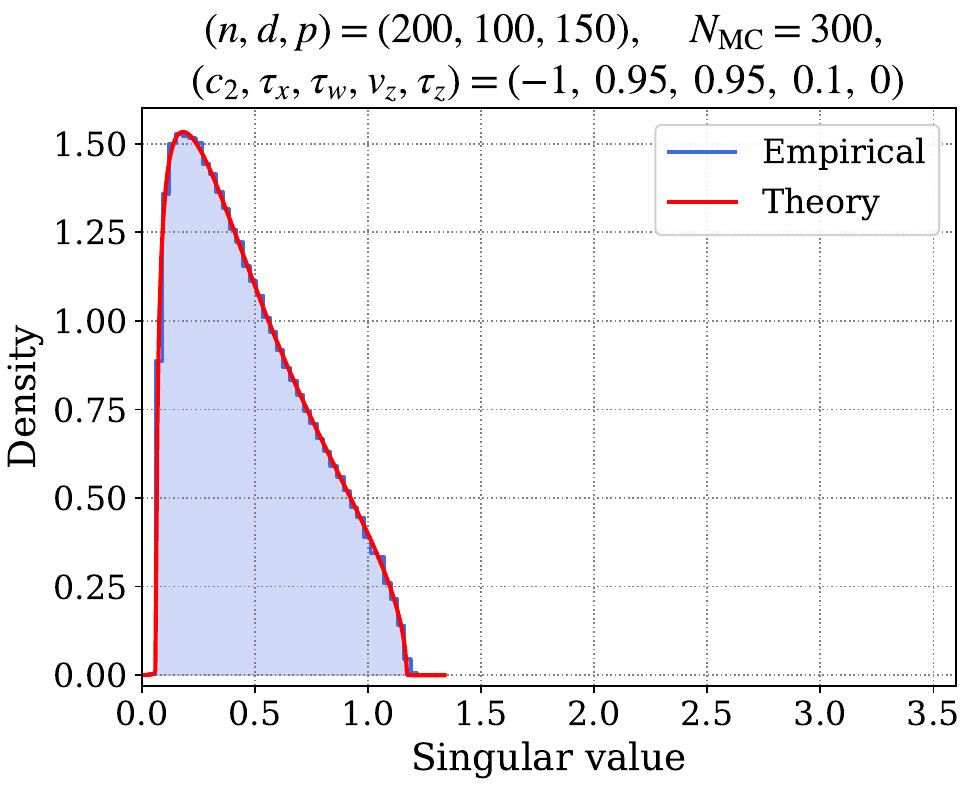}
	\caption{\label{fg:c0czw}Dependence of the spectrum on $\tau_x$ and $\tau_w$ under the condition $T_1^2-4T_2=0$.}
\end{figure}

Figure~\ref{fg:c0czw} shows how the spectrum varies with $\tau_x=\tau_w$ when the remaining parameters are chosen so that the degeneracy condition $T_1^2-4T_2=0$ holds exactly. This setting
corresponds to Case~(c) in Section~\ref{sc:algebraic_reduction}. In contrast to Figure~\ref{fg:ttt}, the spectral support gradually shrinks as $\tau_x$ and $\tau_w$ approach unity, but the density does
not split into two peaks. This behavior agrees with the theoretical picture. The two eigenvalues of $\Sigma_W\Sigma_{CX}$ are degenerate, as shown in \eqref{eq:tT1234r5}. For the parameters used here, with $\tau_x=\tau_w=\tau$, they are given by $\lambda_1=\lambda_2=1-\tau^2$ and decrease together as $\tau$ approaches unity. Consequently, no separation into two distinct spectral scales occurs. At $\tau=1$, both $X$ and $W$ become real, and the two terms $XW$ and $\overline{XW}$ cancel exactly. Thus $Q$ reduces to $Z$ and the eigenvalue distribution of $Q^*Q/n$ is given by the \MP law.

\begin{figure}[tb]
	\centering
	\includegraphics[width=.325\textwidth]{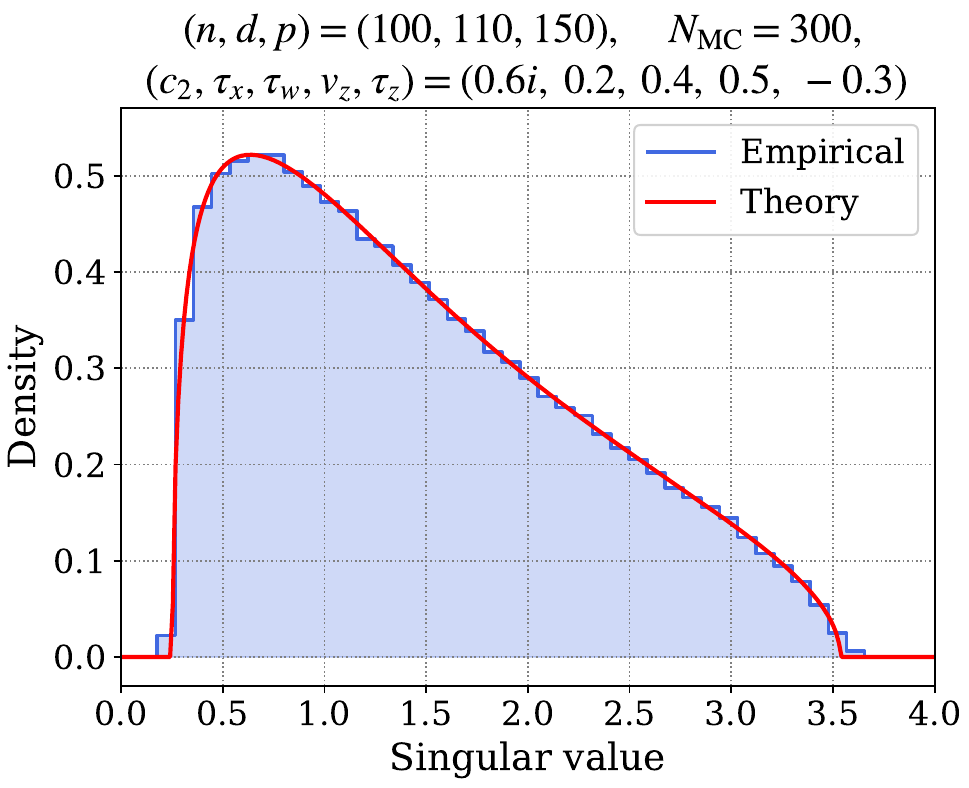}
	\includegraphics[width=.325\textwidth]{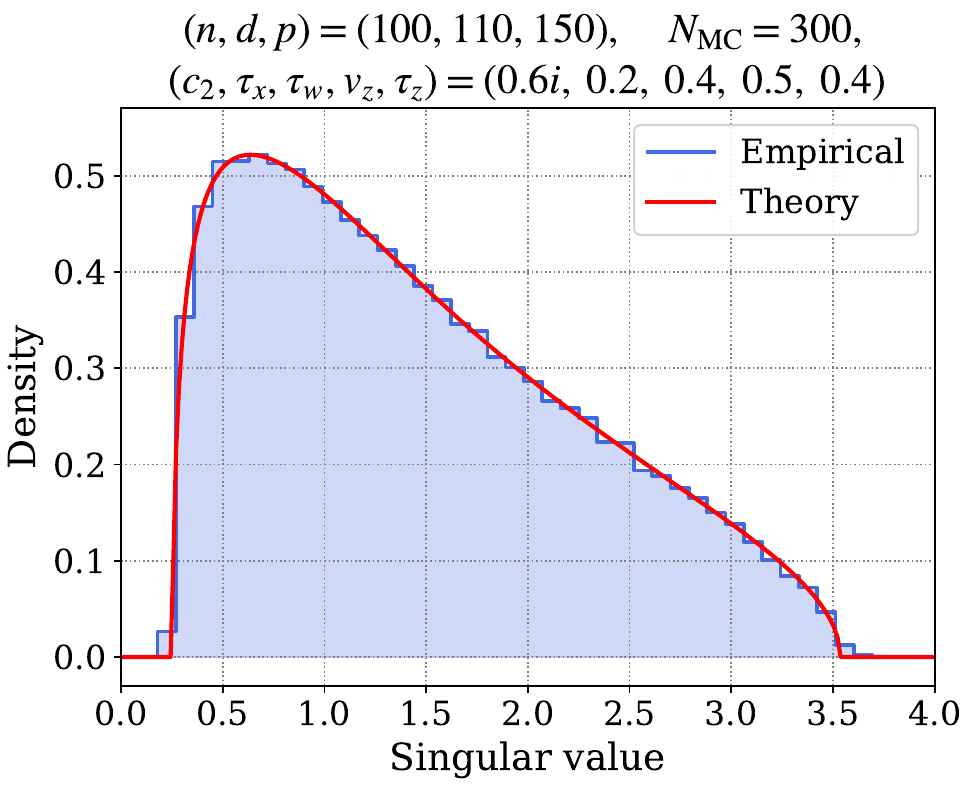}
	\includegraphics[width=.325\textwidth]{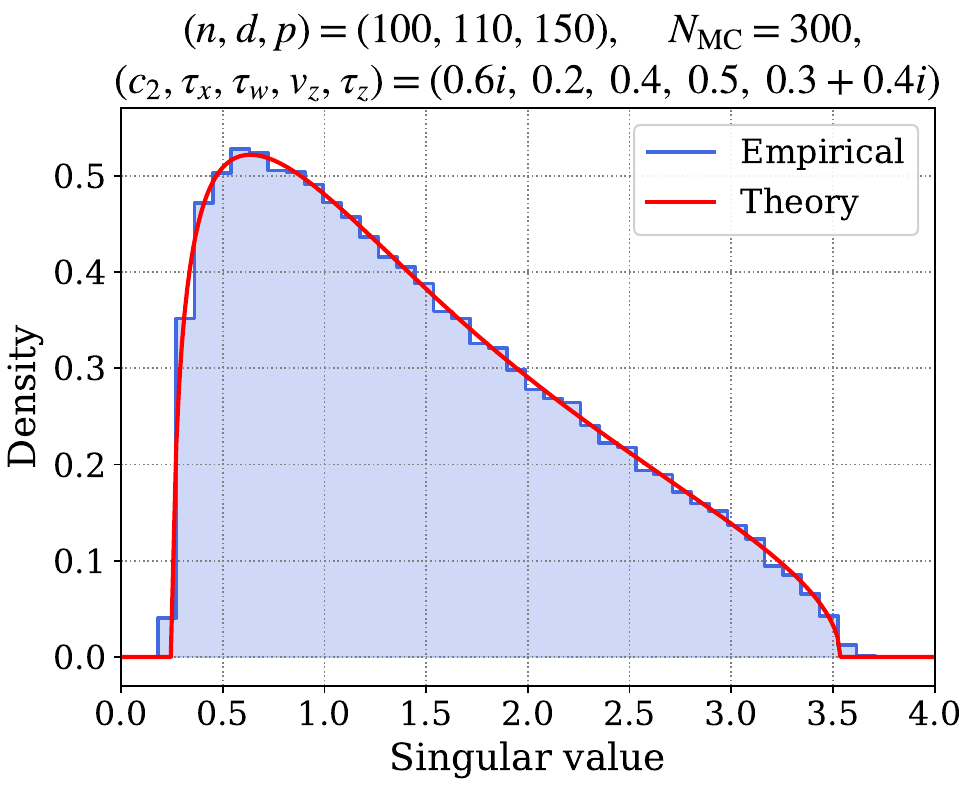}
	\caption{\label{fg:tauzd}Dependence of the spectrum on $\tau_z$.}
\end{figure}

Figure~\ref{fg:tauzd} shows the spectra for varying $\tau_z$ with $v_z=0.5$ fixed. As predicted by the theory, the spectra show no discernible change, confirming that the asymptotic spectral distribution is independent of $\tau_z$.

\begin{figure}[tb]
	\centering
	\includegraphics[width=.325\textwidth]{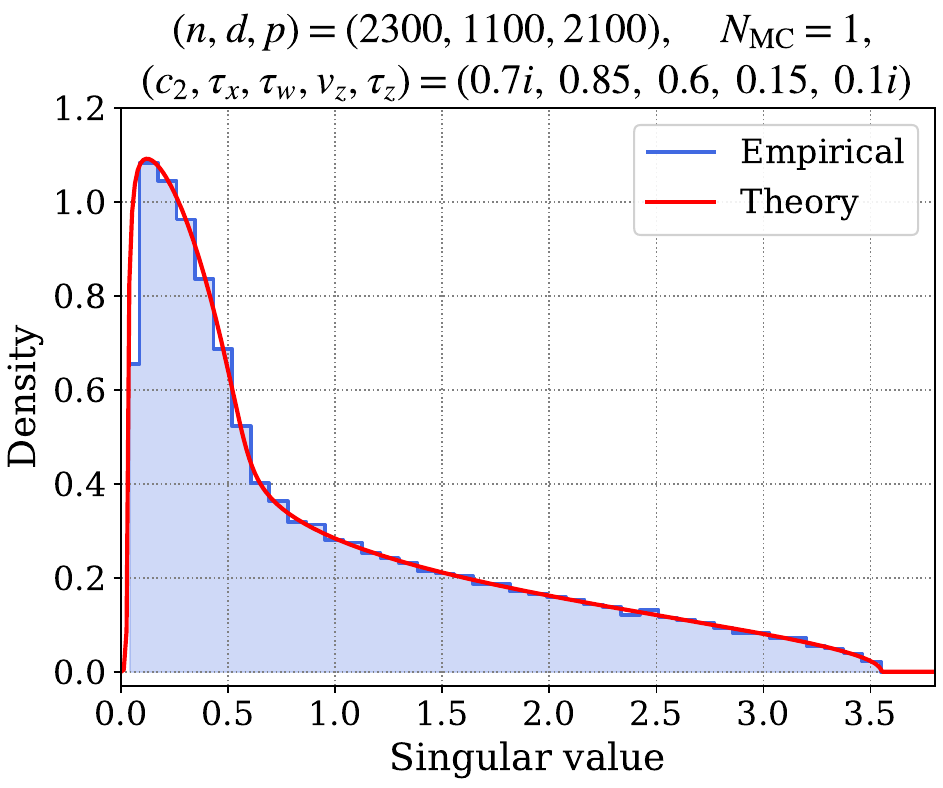}
	\includegraphics[width=.325\textwidth]{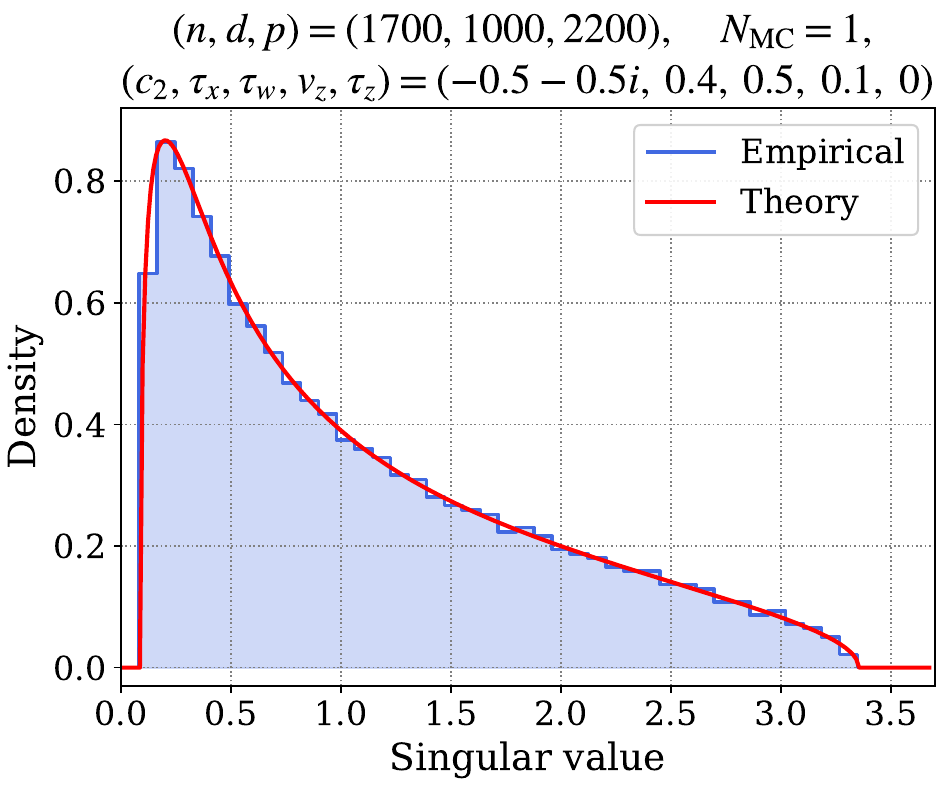}
	\includegraphics[width=.325\textwidth]{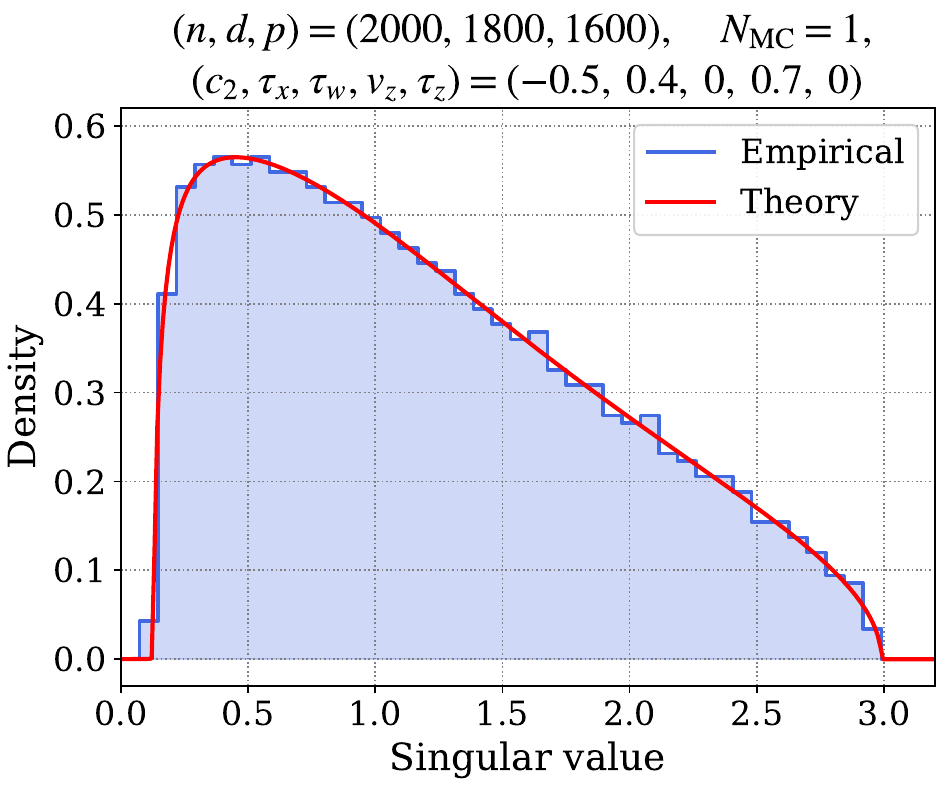}
	\caption{\label{fg:mc1}Empirical spectra from a single realization ($N_{\rm MC}=1$) for three different parameter settings, compared with the corresponding theoretical predictions.}
\end{figure}

Figure~\ref{fg:mc1} compares the theoretical predictions with the empirical spectra obtained from a single realization, without ensemble averaging, for three different choices of the model parameters. We use large matrices whose dimensions are of order $10^3$. Remarkably, the empirical spectra agree closely with the theoretical predictions even for a single realization, providing numerical support for the self-averaging of the empirical spectral distribution.

\begin{figure}[tb]
	\centering
	\includegraphics[width=.325\textwidth]{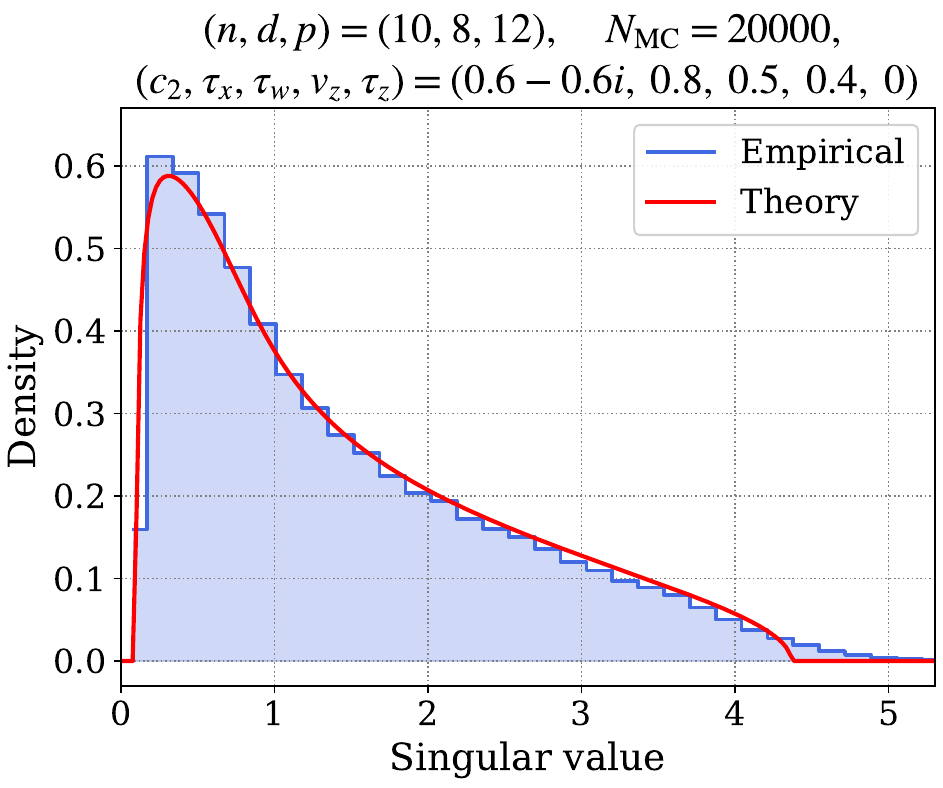}
	\includegraphics[width=.325\textwidth]{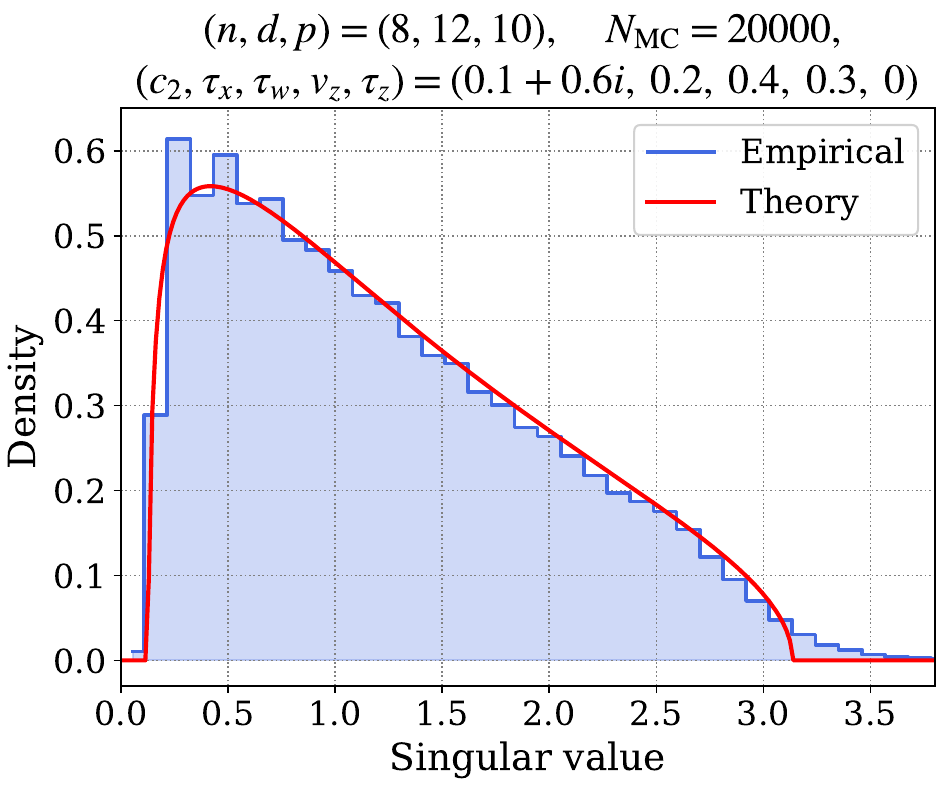}
	\includegraphics[width=.325\textwidth]{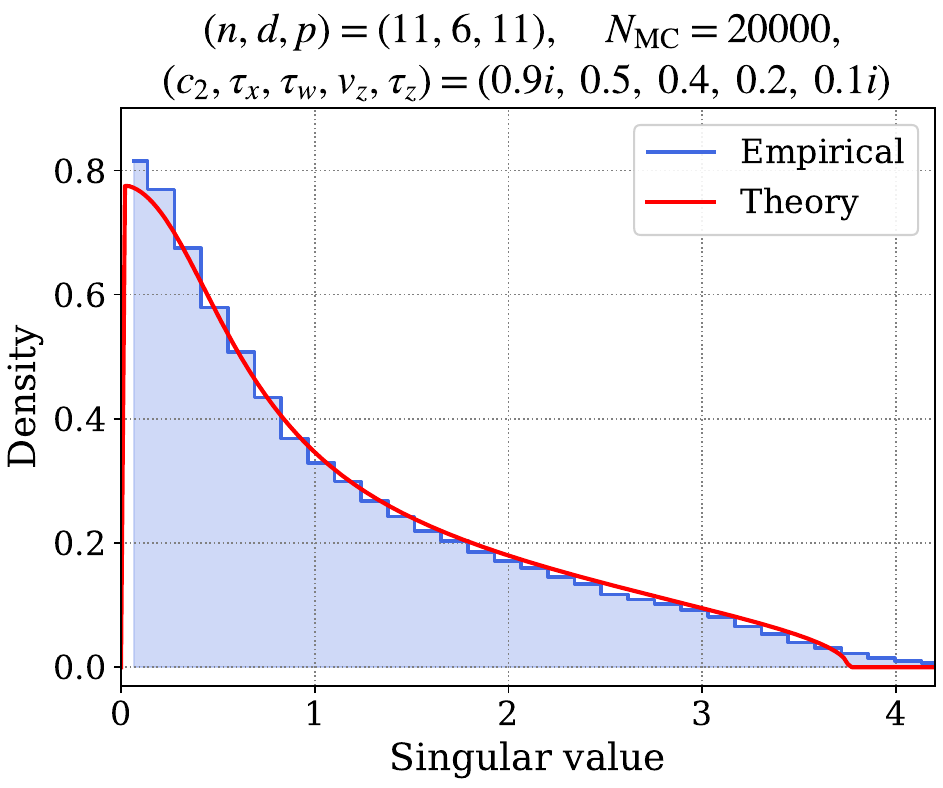}
	\caption{\label{fg:smllr}Empirical spectra for small matrix dimensions of order $10$ under three different parameter settings, compared with the corresponding theoretical predictions.}
\end{figure}

It is interesting to ask how quickly the empirical spectra approach our asymptotic predictions as the matrix dimensions increase. To examine this question, we computed the spectra of substantially smaller matrices for three different parameter settings chosen without fine-tuning, as shown in Figure~\ref{fg:smllr}. Since spectral fluctuations are pronounced at these small dimensions, we used a large number of Monte Carlo realizations, $N_{\rm MC}=20{,}000$, to obtain stable empirical densities. Although a small nonzero tail extends beyond the asymptotic support, the convergence is remarkably fast: even matrix dimensions of order $10$ already show satisfactory agreement with the theoretical predictions.

\begin{figure}[tb]
	\centering
	\includegraphics[width=.325\textwidth]{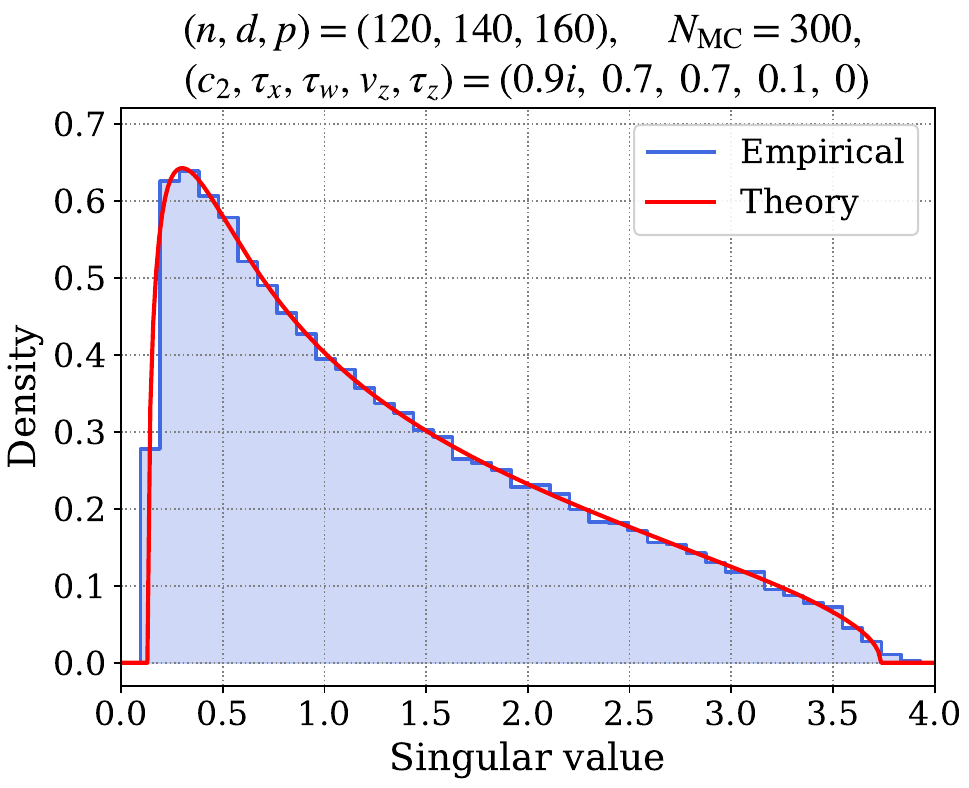}
	\includegraphics[width=.325\textwidth]{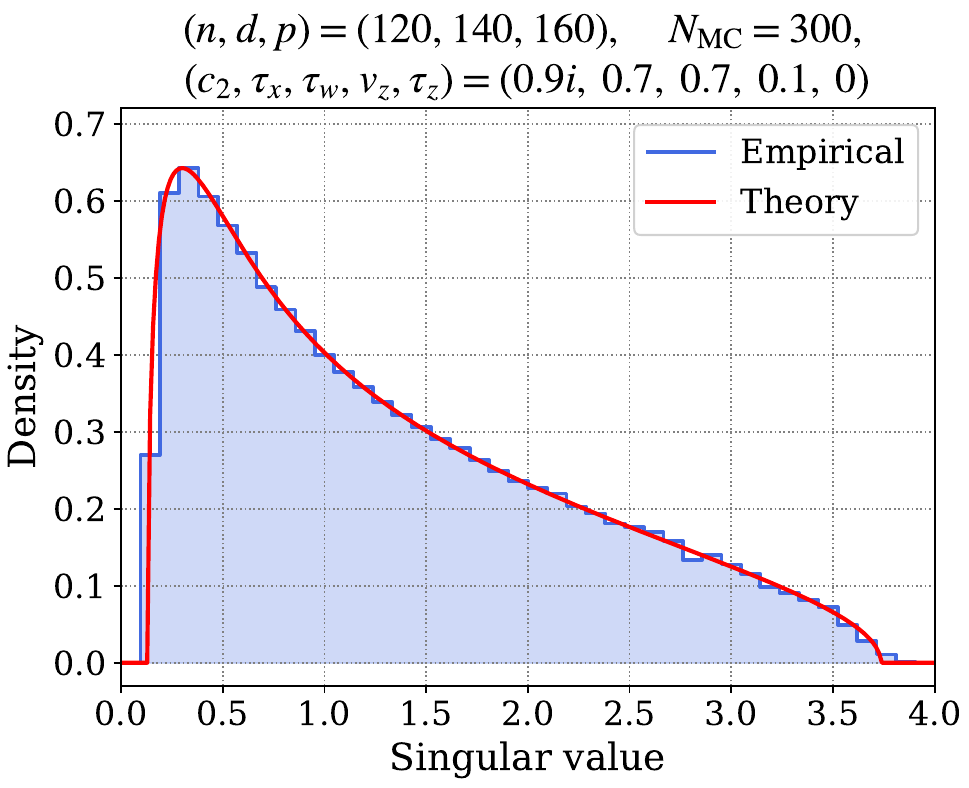}
	\includegraphics[width=.325\textwidth]{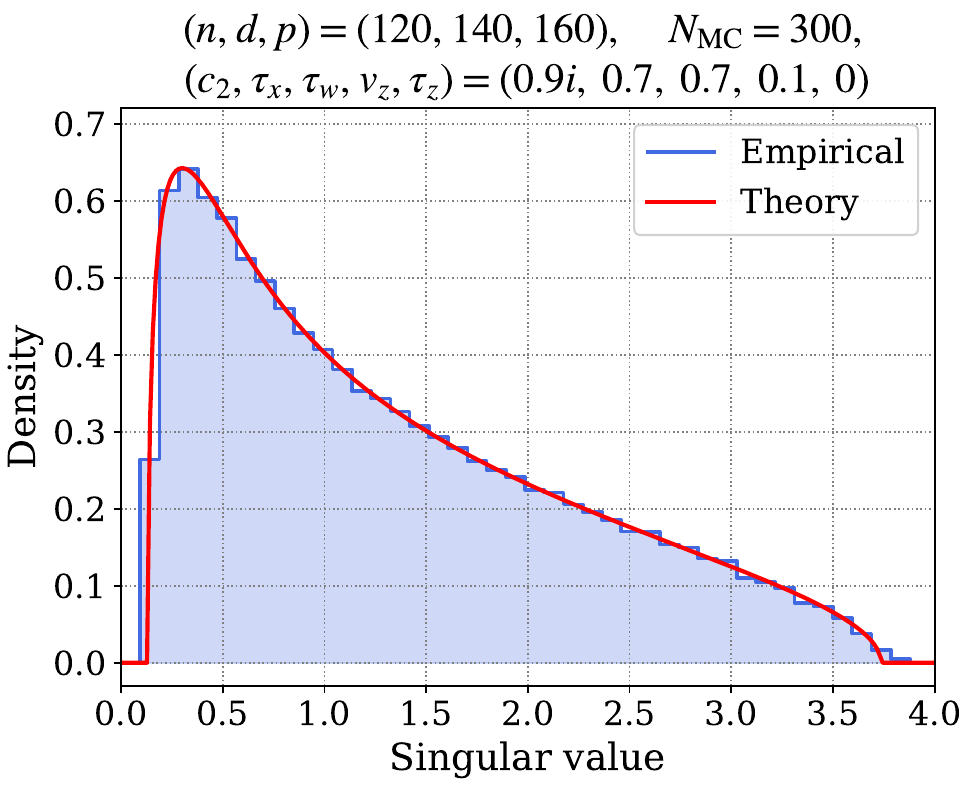}
	\put(-346,60){Gaussian}
	\put(-210,60){Rademacher}
	\put(-60,60){Uniform}
	\vspace{3mm}\\
	\includegraphics[width=.325\textwidth]{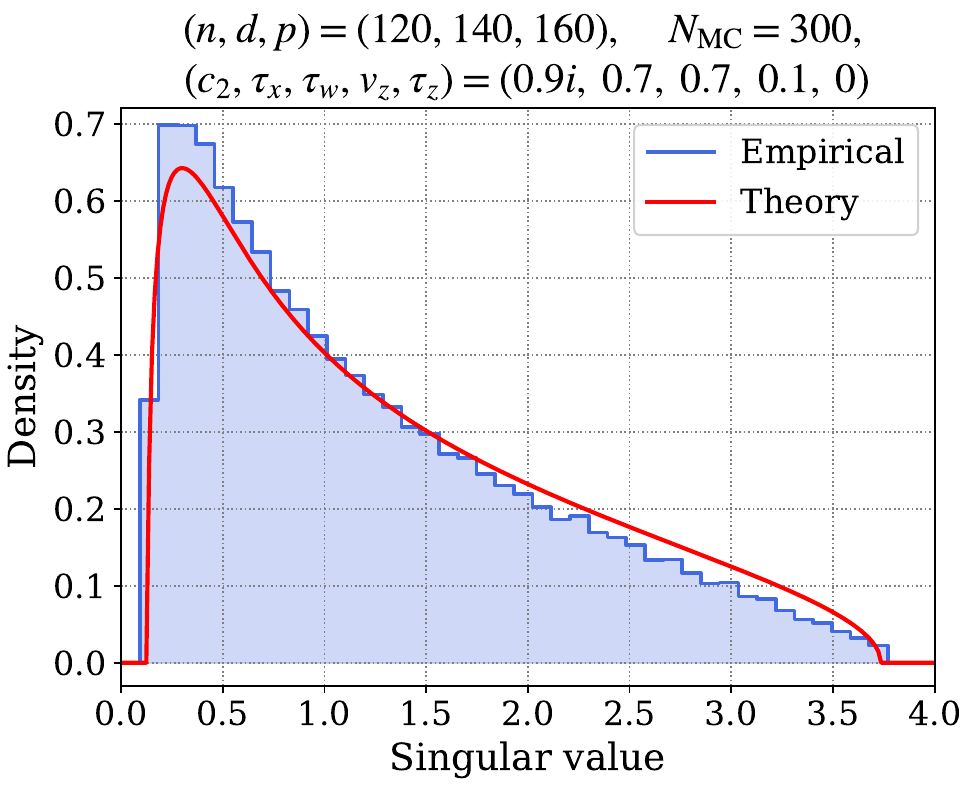}
	\includegraphics[width=.325\textwidth]{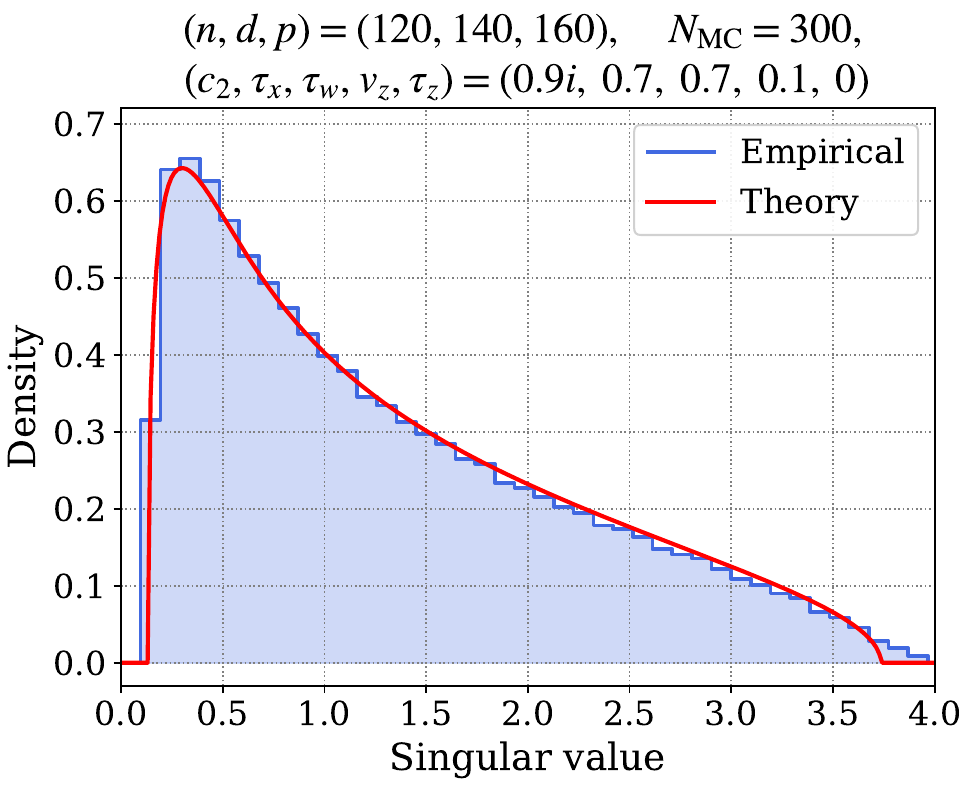}
	\includegraphics[width=.325\textwidth]{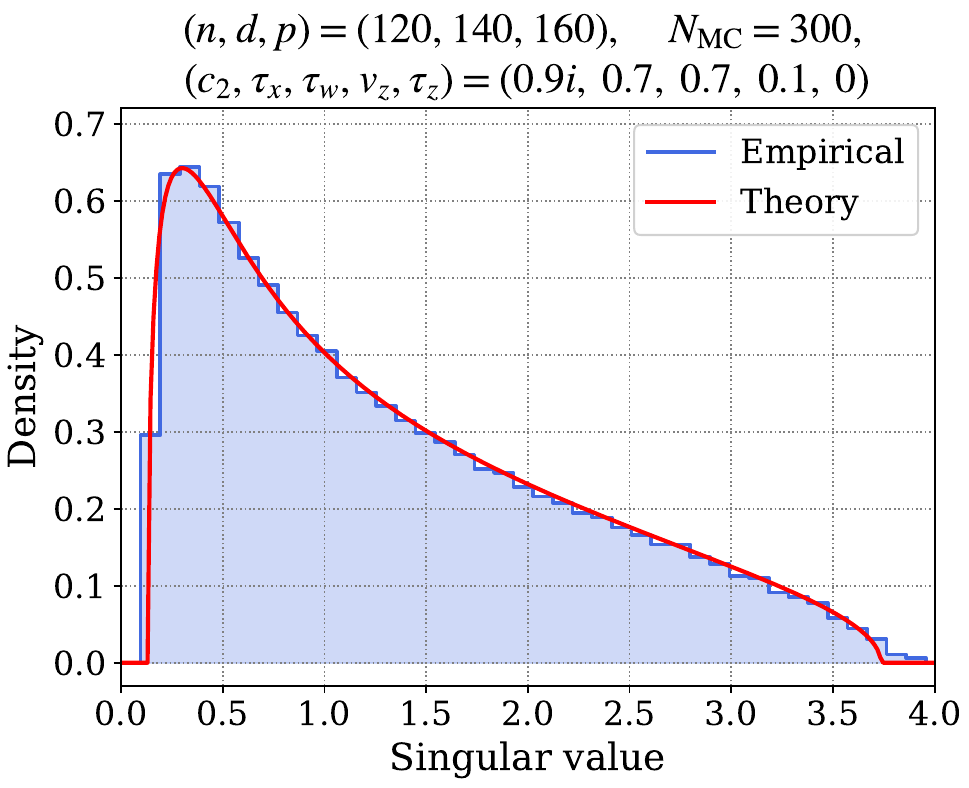}
	\put(-345,58){$\Large\substack{\text{Student-$t$}\\(\nu=3)}$}
	\put(-201,58){$\Large\substack{\text{Student-$t$}\\(\nu=4)}$}
	\put(-58,58){$\Large\substack{\text{Student-$t$}\\(\nu=5)}$}
	\caption{\label{fg:gNg0}Comparison of spectra generated from different entry distributions. The theoretical curve is identical in all panels.}
\end{figure}

\paragraph{Universality beyond Gaussian distributions}
We finally turn to the universality of our results beyond Gaussian entry distributions. All the numerical results presented so far were generated using Gaussian random matrices, in accordance with the ensemble definition in Section~\ref{sc:probset}. However, the derivation of the self-consistent equations in Section~\ref{sc:loo} does not explicitly rely on Gaussianity. More precisely, it relies on the assumption stated in Section~\ref{sc:mainres}, together with the facts that centered linear terms vanish in probability and that quadratic forms concentrate around their corresponding traces. These concentration results are commonly established under suitable moment conditions (see \textcite[Lemma~B.26]{Bai2010book} and \autocite[Lemma~2.11]{Couillet_Liao_2022}), which are satisfied by Gaussian distributions. We have not proved, however, the precise conditions under which our derivation extends to more general entry distributions. 

To examine the universality of the limiting spectral density, we computed the spectra using several different entry distributions for $X$, $W$, and $Z$, while keeping the variance, pseudo-variance, matrix dimensions, and all other model parameters fixed. The results are shown in Figure~\ref{fg:gNg0}. In our numerical implementation, we first draw independent real random variables from a univariate distribution with mean zero and variance one. We then transform pairs of these variables into complex random variables with prescribed variance $v\in\RR_+$ and pseudo-variance $\tau\in\CC$ using Cholesky decomposition of the corresponding $2\times2$ covariance matrix. 

In the top row of Figure~\ref{fg:gNg0}, the Gaussian, Rademacher, and uniform distributions all show excellent agreement with the theoretical prediction. To investigate the effects of heavy tails, we show in the bottom row the results for Student-$t$ distributions with degrees-of-freedom parameter $\nu\in\{3,4,5\}$.%
\footnote{The tail becomes lighter as $\nu$ increases, and the distribution converges to the Gaussian distribution in the limit $\nu\to\infty$.}

Large isolated values occur in all three Student-$t$ cases, reaching as high as $76$ for $\nu=3$. To keep the bulk visible, the bottom-row panels restrict the displayed range; all values are retained in the underlying simulations.

The agreement with the theoretical prediction is excellent for $\nu=4$ and $5$, whereas a clear discrepancy remains for $\nu=3$. This discrepancy is found to gradually diminish as the matrix dimensions increase (not shown). Overall, the experiments suggest that the bulk law extends well beyond the Gaussian setting, although the precise universality conditions and the behavior of heavy-tailed extreme eigenvalues remain open.

\section{Conclusions}
\label{sc:conc}

In this work, we studied the asymptotic spectral distribution of $Q^*Q/n$ for a complex product random matrix with conjugate mixing and additive noise. By constructing a block pencil and applying a leave-one-block-out analysis, we derived a closed system of self-consistent equations for its limiting Stieltjes transform. In the generic case, these equations reduce to a sextic equation, in contrast to the quartic equation arising in the single-channel product-plus-noise model. The higher algebraic degree reflects the additional correlations introduced by the conjugate channel. 

We further classified the regimes in which the sextic equation reduces to lower-degree equations in terms of $T_1$, $T_2$, $T_1^2-4T_2$, and $v_z$. This classification recovers the \MP, product-Wishart, and Pennington--Worah distributions as special cases, while the generic conjugate-mixing regime yields the sextic and quintic spectral equations associated with the two-channel structure. Extensive numerical experiments show excellent agreement with the theoretical predictions, even for small matrix dimensions and single large realizations. They also suggest that the bulk law remains universal beyond Gaussian entry distributions. 

For future work, it would be important to place the present leave-one-out derivation on a fully rigorous foundation by establishing the assumed self-averaging and resolvent stability properties and by identifying precise moment conditions for universality beyond Gaussian distributions. Although free probability offers a possible alternative framework \autocite{Voiculescu1991,Mingo2017book}, its application to non-holomorphic products with improper complex variables is nontrivial; whether it can provide an independent derivation of the present spectral law is an interesting open question. Extensions to products of more than two matrices and more general conjugate-coupled architectures would also be worthwhile. We expect that the spectral theory developed here will be useful for analyzing the performance of complex nonlinear neural networks with non-holomorphic activations.

\section*{Acknowledgments}

The author used GPT-5.6 to assist with English-language editing and to improve the clarity and presentation of the manuscript. The tool was not used as a source of mathematical results. All AI-assisted revisions were reviewed by the author, who takes full responsibility for the final content.

\begin{appendices}
\section{Fundamentals of the Schur Complement}
\label{app:schur}

We briefly review the algebraic structure of the Schur complement \autocite{Zhang2005schur}. Consider an arbitrary square matrix partitioned into a $2 \times 2$ block matrix ${M}$ as follows:
\al{
	{M} = \bmat
	{A} & {B} \\
	{C} & {D}
	\emat.
}
Assuming that the lower-right block ${D}$ is invertible, the matrix
\al{
	{M} / {D} := {A} - {B}{D}^{-1}{C}
	\label{eq:schur_def}
}
is defined as the \emph{Schur complement} of block ${D}$ in ${M}$. Note that the symbol $/$ here denotes the Schur complement operation, not a standard matrix inverse. A salient property of the Schur complement is that its inverse corresponds to the upper-left block of the inverse matrix ${M}^{-1}$. That is,
\al{
	M^{-1} = 
	\bmat
		(M/D)^{-1} & -(M/D)^{-1} B D^{-1} 
		\\
		-D^{-1} C (M/D)^{-1} & 
		D^{-1} + D^{-1} C (M/D)^{-1} B D^{-1} 
	\emat\,.
	\label{eq:schur_inverse}
}

We use the block inversion formula in two ways. First, the \textit{linearization trick} \autocite{Mingo2017book} leverages the block inversion property \eqref{eq:schur_inverse} in reverse. Specifically, it is an approach where one constructs a larger block matrix ${M}$ composed of linear elements in such a way that the nonlinear target matrix of interest naturally emerges in the form of ${A} - {B}{D}^{-1}{C}$. As a toy example, suppose we wish to evaluate the resolvent $G(z) := (-z\1 + {A}^*{A})^{-1}$ for $z \in \CC\setminus\RR_{\geq 0}$. A suitable block pencil is given by
\al{
	{M} = \bmat
		-z\1 & {A}^* \\
		{A} & -\1
	\emat.
}
Applying the Schur complement with respect to the lower-right block $-\1$, we obtain ${M} / (-\1) = -z\1 + {A}^*{A}$. Thus, according to \eqref{eq:schur_inverse}, the upper-left block of the inverse matrix ${M}^{-1}$ is given by $(-z\1 + {A}^*{A})^{-1}$, which is exactly the target resolvent $G(z)$.

Second, as a direct corollary of the Schur complement formula, we can derive the leave-one-out (or leave-one-block-out) resolvent identity, which plays a central role in random matrix theory \autocite{Couillet_Liao_2022}. Consider an arbitrary square matrix ${M}$. By explicitly isolating a specific index $i$ (or an extended block structure, such as a $2 \times 2$ block) and denoting the set of all remaining indices collectively as $\setminus i$, we can rigorously partition ${M}$ into the following $2 \times 2$ macro-block matrix:
\al{
	{M} = \bmat
		{M}_{ii} & 
		{M}_{i, \setminus i} 
		\\
		{M}_{\setminus i, i} & 
		{M}_{\setminus i, \setminus i}
	\emat. 
	\label{eq:loo_m_split}
}
Here, ${M}_{ii}$ is the submatrix corresponding to the isolated element(s), while ${M}_{i, \setminus i}$ and ${M}_{\setminus i, i}$ represent the row and column vectors (or padding blocks) that couple the index $i$ to the rest of the system. The large lower-right block ${M}_{\setminus i, \setminus i}$ represents the submatrix obtained by entirely removing the components associated with index $i$ from the original system. We define the inverse matrix (resolvent) as ${G} = {M}^{-1}$ and partition it conformably with the block structure in \eqref{eq:loo_m_split}:
\al{
	{G} = \bmat
		{G}_{ii} & {G}_{i, \setminus i} 
		\\
		{G}_{\setminus i, i} & 
		{G}_{\setminus i, \setminus i}
	\emat.
}
Assuming that ${M}_{\setminus i, \setminus i}$ is invertible, applying \eqref{eq:schur_inverse} yields 
\al{
	{G}_{ii} & = \big( {M}_{ii} 
	- {M}_{i, \setminus i} {G}^{\s{i}}
	{M}_{\setminus i, i} \big)^{-1},
	\label{eq:loo_identity_final}
}
where we have introduced the resolvent matrix ${G}^{\s{i}} \equiv \mkakko{ {M}_{\setminus i, \setminus i}}^{-1}$. This is the fundamental leave-one-out identity.

\section{Degeneracy Condition for the Sextic Spectral Equation}
\label{app:fac}

This appendix presents a detailed analysis of the degeneracy condition $T_1^2-4T_2=0$ with $T_2>0$ discussed in Section~\ref{sc:algebraic_reduction}. We refer the reader to \eqref{eq:def_T1} and \eqref{eq:def_T2} for the definitions of $T_1$ and $T_2$. Under this condition, the sextic equation factorizes into quadratic and quartic factors, as shown in \eqref{eq:AG43DF9}. 

We begin with the proper limit $\tau_x=\tau_w=0$. In this limit, $T_1=(|c_1|^2+|c_2|^2)v_xv_w$, $T_2=|c_1|^2|c_2|^2 v_x^2v_w^2$, and hence $T_1^2-4T_2 = \mkakko{|c_1|^2-|c_2|^2}^2 v_x^2v_w^2$. The degeneracy condition is therefore satisfied if and only if $|c_1|=|c_2|$. Since $c_1$ and $c_2$ then differ only by a phase, their relative phase can be absorbed by a phase redefinition of $X$ or $W$, so that the widely linear response $c_1XW + c_2\bar{XW}$ takes the form $2\sqrt{c_1c_2}\,\Re(XW)$. Writing $X=X_{\rm R} + i X_{\rm I}$, $W=W_{\rm R} + iW_{\rm I}$, where $X_{\rm R}, X_{\rm I}, W_{\rm R}$ and $W_{\rm I}$ are independent real Gaussian matrices, we obtain
\al{
	2\sqrt{c_1c_2}\,\Re(XW) 
	= 2\sqrt{c_1c_2}\mkakko{X_{\rm R}W_{\rm R} - X_{\rm I}W_{\rm I}}
	= 2\sqrt{c_1c_2} \bmat X_{\rm R} & - X_{\rm I} \emat 
	\bmat W_{\rm R} \\ W_{\rm I} \emat\,.
}
This has the representative form $AB$ reported in Table~\ref{tb:TT}. The variances within each matrix agree, $\Var(X_{\rm R})=\Var(X_{\rm I})$ and $\Var(W_{\rm R})=\Var(W_{\rm I})$, owing to the circular symmetry of $X$ and $W$.

We next turn to the improper case. By definition, 
\al{ 
	T_1 & = (|c_1|^2+|c_2|^2)v_x v_w
	+ 2\Re\mkakko{c_1\bar{c_2}\tau_x \tau_w} 
	\\ 
	& \geq (|c_1|^2+|c_2|^2)v_xv_w 
	- 2|c_1||c_2||\tau_x||\tau_w| \,, 
	\label{eq:47d8QQ17} 
	\\ 
	\sqrt{T_2} &= |c_1||c_2| 
	\sqrt{v_x^2-|\tau_x|^2} \sqrt{v_w^2-|\tau_w|^2}\,. 
} 
The Cauchy--Schwarz inequality gives 
\al{ 
	|\tau_x||\tau_w| + \sqrt{v_x^2-|\tau_x|^2} 
	\sqrt{v_w^2-|\tau_w|^2} \leq v_xv_w\,. 
	\label{eq:cs2eeee}
} 
Combining these relations, we obtain 
\al{
	T_1-2\sqrt{T_2} & \geq 
	(|c_1|^2+|c_2|^2)v_xv_w - 2|c_1||c_2| 
	\mkakko{ |\tau_x||\tau_w| + \sqrt{v_x^2-|\tau_x|^2} 
	\sqrt{v_w^2-|\tau_w|^2} } 
	\\ 
	& \geq (|c_1|^2+|c_2|^2)v_xv_w -2|c_1||c_2|v_xv_w 
	\\ 
	& = (|c_1|-|c_2|)^2 v_xv_w \geq 0\,. 
	\label{eq:Too0oER}
} 
Since $T_1^2-4T_2 = \big(T_1-2\sqrt{T_2}\big)\big(T_1+2\sqrt{T_2}\big)$ and $T_1+2\sqrt{T_2}>0$ when $T_2>0$, the degeneracy condition is equivalent to $T_1-2\sqrt{T_2}=0$. This equality holds if and only if all inequalities in \eqref{eq:47d8QQ17}, \eqref{eq:cs2eeee}, and \eqref{eq:Too0oER} are saturated. Namely, 
\al{ 
	|c_1|=|c_2|>0, \qquad 
	\frac{|\tau_x|}{v_x} = \frac{|\tau_w|}{v_w}, 
	\qquad c_1\bar{c_2}\tau_x\tau_w\in\RR_{\leq 0}\,. 
	\label{eq:cvvno3f}
} 

Assume \eqref{eq:cvvno3f}. Since $|c_1|=|c_2|$, one might be tempted to use the same real-matrix decomposition as in the proper case and derive the expression $2\sqrt{c_1c_2}\mkakko{X_{\rm R} W_{\rm R} - X_{\rm I}W_{\rm I}}$. Under impropriety, however, $X_{\rm R}$ and $X_{\rm I}$, as well as $W_{\rm R}$ and $W_{\rm I}$, are generally correlated and have unequal variances. We must therefore modify this argument by first applying suitable phase rotations. Write $c_1=c_0 \rme^{i\theta_1}$, $c_2=c_0\rme^{i\theta_2}$, $\tau_x=|\tau_x|\rme^{2i\theta_x}$, and $\tau_w=|\tau_w|\rme^{2i\theta_w}$ with $c_0>0$, and define $X'\equiv\rme^{-i\theta_x}X$ and $W'\equiv\rme^{-i\theta_w}W$. Then their entrywise pseudo-variance reads as
\al{
	\EE\big[X'^2_{i\mu} \big] & = \rme^{-2i\theta_x}
	\EE\big[X_{i\mu}^2 \big]
	= \rme^{-2i\theta_x} \tau_x = |\tau_x|\,,
	\\
	\EE\big[W'^2_{\mu\alpha} \big]
	& = \rme^{-2i\theta_w}\EE\big[W_{\mu\alpha}^2 \big]
	= \rme^{-2i\theta_w} \tau_w = |\tau_w|\,.
}
Since these pseudo-variances are real, the real and imaginary parts $X'_{\rm R}$ and $X'_{\rm I}$ are uncorrelated and hence independent by Gaussianity. The same holds for $W'_{\rm R}$ and $W'_{\rm I}$. Furthermore,
\al{
	c_1 XW + c_2 \bar{XW}
	& = c_0 \rme^{i(\theta_1+\theta_x+\theta_w)} X'W'
	+ c_0 \rme^{i(\theta_2-\theta_x-\theta_w)} \bar{X'W'}
	\\
	& = c_0 \rme^{i(\theta_1+\theta_x+\theta_w)}
	\big(
		X'W'+ \rme^{i(-\theta_1+\theta_2
		-2\theta_x-2\theta_w)}\bar{X'W'}
	\big)\,.
}
Now we use $\rme^{i(-\theta_1+\theta_2 - 2\theta_x-2\theta_w)}=-1$, which follows from the last condition in \eqref{eq:cvvno3f}, to obtain
\al{
	c_1 XW + c_2 \bar{XW}
	& = c_0 \rme^{i(\theta_1+\theta_x+\theta_w)}
	(X'W' - \bar{X'W'})
	\\
	& = 2 i c_0 \rme^{i(\theta_1+\theta_x+\theta_w)}
	\big( X'_{\rm I} W'_{\rm R} + X_{\rm R}'W_{\rm I}' \big)
	\propto \bmat X'_{\rm I} & t^{-1} X_{\rm R}' \emat
	\bmat W'_{\rm R} \\ t W_{\rm I}' \emat 
	\equiv AB\,,
	\label{eq:ABdef}
}
where $t>0$ is an arbitrary constant. We next examine the variances:
\al{
	\Var\big(X'_{\rm R}\big) = \frac{v_x+|\tau_x|}{2}\,, & \quad
	\Var\big(X'_{\rm I}\big) = \frac{v_x-|\tau_x|}{2}\,,
	\\
	\Var\big(W'_{\rm R}\big) = \frac{v_w+|\tau_w|}{2}\,, & \quad
	\Var\big(W'_{\rm I}\big) = \frac{v_w-|\tau_w|}{2}\,.
}
Then $\Var\big(X'_{\rm I}\big)=\Var\big(t^{-1}X'_{\rm R}\big)$ and $\Var\big(W'_{\rm R}\big)=\Var\big(t W'_{\rm I}\big)$ can be achieved simultaneously by choosing
\al{
	t = \sqrt{\frac{1 + |\tau_x|/v_x}{1-|\tau_x|/v_x}}
	= \sqrt{\frac{1 + |\tau_w|/v_w}{1-|\tau_w|/v_w}} \,.
}
Since $T_2>0$ is assumed, the denominators are positive. The second equality follows from \eqref{eq:cvvno3f}. With this value of $t$, $A$ and $B$ are independent real random matrices, each with \mbox{i.i.d.} entries. We thus conclude that, as reported in Table~\ref{tb:TT} for Cases~(c) and~(d), the two-channel model $c_1 XW + c_2 \bar{XW}$ reduces to a single real product with inner dimension $2d$ and the aspect ratios $\phi/2$ and $\psi/2$. This explains algebraically why the equations in Cases~(c) and~(d) agree with those in Cases~(e) and~(f), respectively, after the corresponding substitution.
\end{appendices}

\begin{refcontext}[sorting=nyt]
\printbibliography
\end{refcontext}
\end{document}